\documentclass{scrartcl}

\usepackage{verbatim} 
\usepackage[utf8]{inputenc}
\usepackage{mathbbol,mathrsfs,amsfonts,amssymb,amsthm,stmaryrd}
\usepackage{amsmath}
\usepackage[all]{xypic}
\usepackage{mathpazo}
\usepackage[mathcal]{euscript}
\usepackage{graphics,supertabular,color}
\usepackage[colorlinks,hyperindex]{hyperref}
\hypersetup{linkcolor=blue, citecolor=blue}

\newcommand{\C}{\mathbb{C}}
\newcommand{\R}{\mathbb{R}}

\newcommand{\Z}{\mathbb{Z}}
\newcommand{\U}{{\rm U}}

\newtheoremstyle{ansgarthmstyle}{}{}{\it}{}{\bf}{}{0.4cm}{}

\theoremstyle{ansgarthmstyle}

\newtheorem{exa}{Example}[section]

\newtheorem{defi}{Definition}[section]
\newtheorem{prop}{Proposition}[section]
\newtheorem{rem}{Remark}[section]
\newtheorem{lem}{Lemma}[section]
\newtheorem{thm}{Theorem}[section]
\newtheorem{cor}{Corollary}[section]

\newenvironment{pf}{\begin{proof}[\textrm{{\bf Proof\ :\!}}]} {\end{proof}}

\begin{document}

\title{The Global Topology of Pontrjagin Duality}

\author{Ansgar Schneider
\footnote{Max-Planck-Institut f\"ur Mathematik, Vivatsgasse 7, 53111 Bonn, ansgar@mpim-bonn.mpg.de} 
}

\maketitle

\begin{abstract}
\noindent
Pontrjagin duality is implemented
in the framework of fibre bundles.
By means of Pontrjagin duality triples a Fourier 
transform is defined by a pull-push construction operating 
on sections of line bundles. This yields an 
isomorphism of Hilbert $C^*$-modules
which generalises the classical isomorphism 
between the group $C^*$-algebra of a group
and the continuous functions vanishing at infinity 
on the dual group. 
\end{abstract}

\tableofcontents

\section{{Introduction}}
\noindent
The main tool of classical harmonic analysis on abelian groups 
is the Fourier transform. It maps  a function $\alpha:G\to \C$ 
on a  locally compact abelian group $G$ 
to  a function $\hat \alpha:\widehat G\to \C$
on the dual group $\widehat G:={\rm Hom}(G,\U(1))$
of  $G$. Explicitely, it is given by the formula
\begin{eqnarray*}
\hat\alpha(\chi):=\int_G \alpha(g){\langle g,\chi\rangle}\ dg, \quad\chi \in \widehat G,
\end{eqnarray*}
where $\langle\ ,\ \rangle\colon G\times \widehat G\to \U(1)$ is the 
canonical pairing and $dg$ is the Haar measure on $G$.
As it turns out, the Fourier transform of an integrable
function is a continuous function vanishing at infinity, 
and in fact it extends to an isomorphism of $C^*$-algebras
\begin{eqnarray}\label{EqOfAlgebrasinnerIntro}
\hat{\ }:C^*(G)&\cong&C_0(\widehat G)
\end{eqnarray}
between the group $C^*$-algebra of $G$ (which 
contains the integrable functions as a dense subspace)
and the continuous functions on $\widehat G$ which vanish at infinity.

We  re\"interpret these data in a bundle theoretic set-up.
By addition the space $G$ has a free (right) action of the group 
$G$  which means that the quotient map of this action
$G\to\ast$ is a trivial $G$-principal fibre bundle 
over the one-point space $\ast$. 
(Recall that a $G$-principal 
fibre bundle $E\to B$ over a topological space $B$
is a space $E$ with a free (right) $G$-action and a
homoeomorphism $E/G\cong B$ such that the quotient 
map $E\to B$ has local sections.) 
A complex valued function $\alpha: G\to \C$ corresponds 
canonically to a section 
$$
\resizebox*{!}{2,6cm}{$
\xymatrix{
L\ar[d]\\
G\ar@/_0.5cm/[u]_\alpha\ar[d]\\
\ast
}$
}
$$
of the  line bundle $L:=(G\times\U(1))\times_{\U(1)}\C$
associated to the trivial $\U(1)$-principal fibre bundle $G\times \U(1)\to G$.
(Recall that the associated line bundle of 
a $\U(1)$-principal  bundle $F\to E$ is just the quotient 
of $F\times \C$ by the diagonal $\U(1)$-action:
$((x,z), t)\mapsto (x\cdot t, t^{-1}z)$.)
Thus, on the domain side of the Fourier transform
we have an underlying topological 
object $G\times \U(1)\to G\to \ast$ which
we call the trivial pair over the one-point space.
On the range side of the Fourier transform 
we have by analogy a trivial dual pair which is the sequence
$\widehat G\times \U(1)\to\widehat G\to \ast$.
Together these two objects form a diagram 
of principal fibre bundles over the one-point space
\begin{equation*}
\resizebox*{!}{2,5cm}{
\xymatrix{
&\times\ar[rd]&\hspace{-2,2cm}G\ \ \ \ \U(1)&
&&\times\ar[ld]&\hspace{-2cm}\widehat G\ \ \ \ \U(1)\\
&&G\ar[rd]&&\widehat G\ar[ld]&&\\
&&&\ast&&&
}}.
\end{equation*}
We extend this diagram step by step to
a diagram of pullbacks
\begin{equation*}
\resizebox*{!}{4cm}{
\xymatrix{
&G\times\widehat G\times \U(1)\ar[rd]\ar[ld]&
&G\times\widehat G \times\U(1)\ar[rd]\ar[ld]&\\
G\times\U(1)\ar[rd]&&G\times \widehat G\ar[rd]\ar[ld]&
&\widehat G\times\U(1)\ar[ld] \\
&G\ar[rd]&&\widehat G\ar[ld]&\\
&&\ast&&
}}.
\end{equation*}
The canonical pairing
$\langle\ ,\ \rangle\colon G\times\widehat G\to \U(1)$
gives rise to a $\U(1)$-principal fibre bundle morphism 
\begin{eqnarray*}
\xymatrix{
G\times\widehat G\times\U(1)\ar[dr]\ar[rr]^\pi&&
G\times\widehat G\times\U(1)\ar[dl]\\
&G\times\widehat G
}
\end{eqnarray*}
just given by $\pi(g,\chi,t):=(g,\chi, \langle g,\chi\rangle t).$
Inserting this isomorphism into the previous 
diagram we obtain a diagram
\begin{equation}\label{DiagDatTriviInnerIntro}
\resizebox*{!}{4cm}{
\xymatrix{
&G\times\widehat G\times \U(1)\ar[rd]\ar[ld]\ar[rr]^{\pi}&
&G\times\widehat G \times\U(1)\ar[rd]\ar[ld]&\\
G\times\U(1)\ar[rd]&&G\times \widehat G\ar[rd]\ar[ld]&
&\widehat G\times\U(1)\ar[ld] \\
&G\ar[rd]&&\widehat G\ar[ld]&\\
&&\ast&&
}}
\end{equation}
which we call the trivial Pontrjagin duality triple
over the one-point space.
By this diagram we can define a Fourier transform 
in disguise which maps the sections of 
the associated line bundle $L\to G$
to the sections of the 
``dual" 
line bundle 
$\widehat L:= (\widehat G\times\U(1))\times_{\U(1)}\C\to \widehat G$.
Explicitly, the Fourier transform of a section 
$\alpha:G\to  L$ 
is the section $\hat\alpha:\widehat G\to \widehat L$
given by
\begin{eqnarray*}\label{EqDieFTistimmergut}
\hat \alpha(\chi):=\int_G {\rm pr}_{\widehat L}\Big(
\pi^\C\big(\alpha(h),\chi\big) \Big)\ dh,
\end{eqnarray*}
where $\pi^\C:L\times \widehat G\to G\times\widehat L$ 
is the isomorphism of line bundles which is induced by $\pi$,
and ${\rm pr}_{\widehat L}:G\times \widehat L\to \widehat L$ is the 
projection.
By renaming the objects in (\ref{EqOfAlgebrasinnerIntro}) 
we have an isomorphism  
\begin{eqnarray*}
\hat{\ }:C^*(G,G\times\U(1))\cong\Gamma_0(\widehat G,\widehat L),
\end{eqnarray*}
where  $C^*(G,G\times \U(1))$ is the $C^*$-algebra 
of the trivial pair $G\times\U(1)\to G\to \ast$, and 
$\Gamma_0(\widehat G, \widehat L)$
are the sections vanishing at infinity.

Now, the stage is set for topology. We want to consider 
family versions of diagram (\ref{DiagDatTriviInnerIntro})
glued by the topology of a space $B$ which leads 
us to the notion of a general Pontrjagin duality
triple. A Pontrjagin duality triple is a commutative diagram 
of principal fibre bundles 
\begin{eqnarray}\label{DiagPDTInnerIntro}
	\xymatrix{
	&F\times_B \widehat E\ar[rd]\ar[ld]\ar[rr]^\kappa&
	&E\times_B \widehat F\ar[rd]\ar[ld]&&\\
	F\ar[rd]&\hspace{-1,8cm}\rotatebox {90}{$\circlearrowleft$}\ { \U(1)}
	&E\times_B \widehat E\ar[rd]\ar[ld]&
	&\widehat F\ar[ld]&\hspace{-0.9cm}\rotatebox {90}{$\circlearrowleft$}\ {\U(1)}\\
	&E\ar[rd]&\hspace{-3cm}\rotatebox {90}{$\circlearrowleft$}\ { G}
	&\widehat E\ar[ld] &\hspace{-2,1cm}\rotatebox {90}{$\circlearrowleft$}\ \widehat G&\\
	&&B&&&
	}
\end{eqnarray}
such that the restriction of (\ref{DiagPDTInnerIntro}) 
to each point $b\in B$ is isomorphic to (\ref{DiagDatTriviInnerIntro}). 
Let  $ F^\C:=  F\times_{\U(1)}\C$ and 
$\widehat F^\C:= \widehat F\times_{\U(1)}\C$ 
denote the associated line bundles.
For a horizontally integrable section $\gamma: E\to F^\C$
we can define the Fourier transform based on diagram
(\ref{DiagPDTInnerIntro})  to be a section 
$\hat \gamma:\widehat E\to \widehat F^\C$.
Namely, for $\hat e\in\widehat E$ over $b\in B$ take any 
$e\in E$ also over $b$ and define 
\begin{eqnarray}\label{MhNamenloseGl}
	\hat \gamma(\hat e):=\int_G  {\rm pr}_{\widehat F^\C}\Big(
	\kappa^\C\big(\gamma(e\cdot h),\hat e\big)\Big)\ dh,
\end{eqnarray}
where $\kappa^\C: F^\C\times_B \widehat E\to E\times_B\widehat F^\C$
is the induced isomorphism of the top isomorphism $\kappa$ in
(\ref{DiagPDTInnerIntro}), and ${\rm pr}_{\widehat F^\C}:
E\times_B\widehat F^\C\to \widehat F^\C$ is the 
projection.

We wish to understand the Fourier transform (\ref{MhNamenloseGl})
in the correct $C^*$-algebraic context.
The answer will be that it defines an isomorphism of Hilbert 
$C^*$-modules (Theorem \ref{ThmTheMainThm}).
In fact, these Hilbert $C^*$-modules are Hilbert 
$C^*$-modules of $\U(1)$-equivariant (self) 
Morita equivalences, namely $F$ and $\widehat F$,
between $\U(1)$-central extensions of groupoids.
The task is to construct these groupoids. 
\\

Let us give diagrams (\ref{DiagDatTriviInnerIntro}) and 
(\ref{DiagPDTInnerIntro})  another purely topological 
interpretation.
The dual of the group $G$, understood as a $G$-principal
fibre bundle $G\to\ast$,  is the dual group $\widehat G$,
understood as a $\widehat G$-principal fibre bundle
$\widehat G\to \ast$.
Pontrjagin duality triples can give answer to the 
question of the dual of a general 
$G$-principal fibre bundle $E\to B$.
Note that the a $G$-bundle $E\to B$
is the same amount of data as the 
bundle together with a trivial $\U(1)$-bundle
$E\times\U(1)\to E$ on its total space. 

We introduce some terminology. A pair $F\to E\to B$ is a sequence 
of principal fibre bundles which is locally 
isomorphic to the trivial pair $B\times G\times \U(1)\to B\times G\to B$.
A dual pair $\widehat F\to \widehat E\to B$ is
a sequence of principal fibre bundles 
which is locally isomorphic to the dual trivial pair 
$ B\times\widehat G\times \U(1)\to B\times\widehat G\to B$.
We call $\widehat F\to \widehat E\to B$ a dual of 
$F\to  E\to B$ if one can extend these two to a 
Pontrjagin duality triple (\ref{DiagPDTInnerIntro}).

So, if  $E\to B$ is  a $G$-principal fibre bundle,
then $E\times \U(1)\to E\to B$ is a pair, and 
we are concerned  with the questions of existence 
and uniqueness of duals in the above sense. 

These questions can be answered by analysing 
the topology of the classifying spaces of 
the automorphism groups of the corresponding 
local models.
The automorphism group of the trivial pair over the point 
$G\times \U(1)\to G\to \ast$ is the semi-direct product 
${\rm A_{\rm Par}}:=G\ltimes C(G,\U(1))$, and the 
category of pairs is equivalent to the category 
of ${\rm A_{\rm Par}}$-principal fibre bundles
(Proposition \ref{PropEqOfCatsPairsAndBunAFull}).
The automorphism group of diagram (\ref{DiagDatTriviInnerIntro})
is the semi-direct product ${\rm A_{\rm Pon}}:=G\ltimes(\U(1)\times\widehat G)$,
and the category of Pontrjagin duality triples is 
equivalent to the category of ${\rm A_{\rm Pon}}$-principal fibre bundles
(Proposition \ref{PropEqOfCatsPonAndAutoOfPon}). 
${\rm A_{\rm Pon}}$ is the subgroup
of ${\rm A_{\rm Par}}$ consisting of those $(g,f)\in {\rm A_{\rm Par}}$
which admit a $t\in\U(1)$ and a $\chi\in\widehat G$
such that $f(h)=t\langle h,\chi\rangle$.
The inclusion 
$
{\rm A_{\rm Pon}}\hookrightarrow{\rm A_{\rm Par}}
$
induces a map between the classifying spaces
$
{\rm BA_{\rm Pon}}\to {\rm BA_{\rm Par}}.
$

Expressed in homotopy theoretic terms, 
the question of the existence of a dual of 
a pair $F\to E\to B$ is the question
of the existence of a (homotopy) lift  
of the classifying map $B\to {\rm BA_{\rm Par}}$ 
of $F\to E\to B$: 
\begin{eqnarray*}
\xymatrix{
&&{\rm BA_{\rm Pon}}\ar[d]\\
B\ar[rr]\ar@{..>}[rru]^{\exists ?}&& {\rm BA_{\rm Par}}
}.
\end{eqnarray*}
If such a lift exists, then the question of uniqueness 
is the question whether the (homotopy class of this) 
lift is unique. 
As the topology of the classifying spaces varies with
the group $G$, the answers to these questions depend 
on the group $G$ and are quite different for different $G$
(s. section \ref{SecExtensions}).

\begin{rem}
	The notion of Pontrjagin duality triples is similar to what 
	has been introduced  in \cite{BRS} under the name 
	T-duality triples (see also \cite{Sch1,Sch2,BSST}). 
	The analysis of  Pontrjagin duality triples in this work 
	consists of similar steps as the analysis of the 
	$C^*$-algebraic content of T-duality triples \cite{Sch1,Sch2}, 
	but the necessary tools for their investigation stay 
	on a much more explicit level.
\end{rem}
\vspace{1cm}
\noindent
{\bf Conventions.}
By a {\bf space} $B$ we will always mean a Hausdorff, paracompact topological space. Frequently, we use the word {\bf bundle} 
as an abbreviation for principal fibre bundle.
In  the whole of this work we denote by $G$ a Hausdorff, second countable,
locally compact abelian group. Its dual group is 
$\widehat G:={\rm Hom}(G,\U(1))$. With the compact-open topology
$\widehat G$ is again a Hausdorff, second countable, locally compact
abelian group. We denote by 
$
\langle . , .\rangle : G\times \widehat G\to \U(1)
$
the canonical paring.
The integral of a function $f$ on $G$ with the Haar measure 
is simply denoted by $\int_G f(g)\ dg$.

We refer to \cite{Ru} for the classical theory
of Fourier analysis and to \cite{La} for the 
language of Hilbert $C^*$-modules.
An introduction to the language of stacks
is found in \cite{Hei} or \cite{Moe}.

\vspace{1cm}
\noindent
{\bf Acknowledgement.}
The author is grateful to Ulrich Bunke and Ilya Shapiro for helpful discussions and suggestions.

\section{Pairs}\label{SecPairs}
\noindent
We start with the definition of a pair which is a collection of
$$
\xymatrix{
G\times \U(1)\ar[d]\\
G\ar[d]\\ 
\ast
}
$$
glued together by the topology of a space $B$.

\begin{defi} 
	A {\bf pair} over a space $B$ is a sequence $F\to E \to B$
	of a $G$-bundle $E\to B$ and a $\U(1)$-bundle $F\to E$ 
	subject to the following local triviality axiom:
	There is an open cover $\{U_i\}$ of $B$ together with bundle
	isomorphisms
	\begin{equation}\label{DiagChartOfAPair}
	\xymatrix{
	F|_{E|_{U_i}}\ar[r]\ar[d]& U_i\times G\times \U(1)\ar[d]\\
	E|_{U_i}\ar[r]\ar[d]&U_i\times G\ar[d]\\ 
	U_i\ar[r]^=&U_i
	},
	\end{equation}
	i.e. we require the $\U(1)$-bundle $F$ to be trivialisable over
	the fibres of $E$.
	
	The notion of a {\bf dual pair}  $\widehat F\to\widehat E \to B$
	is defined by the same means, but with $G$ replaced 
	by its dual group $\widehat G$.

	The open cover $\{U_i\}$ of $B$ together with the diagrams 
	(\ref{DiagChartOfAPair}) is called an {\bf atlas} for
	the pair, and each single diagram (\ref{DiagChartOfAPair})
	is called a {\bf chart}.
\end{defi}

We define the category of pairs over 
a space $B$ in the obvious way, i.e.
a morphism of pairs is a diagramm 
$$
\xymatrix{
F\ar[r]\ar[d]& F'\ar[d]\\
E\ar[r]\ar[d]&E'\ar[d]\\ 
B\ar[r]^=&B
}
$$
with the horizontal arrows being  
bundle morphisms. This category clearly is
a groupoid.

Let $F\to E\to B$ be a pair and let $\{ U_i\}$ be an atlas.
By the usual arguments, we obtain transition functions 
on twofold intersections $U_{ij}:= U_i\cap U_j$ 
\begin{eqnarray}\label{EqTransisForPairs}
g_{ji}:U_{ji}\to G, \qquad \zeta_{ji}:U_{ij}\times G\to \U(1)
\end{eqnarray}
which satisfy 
\begin{eqnarray}\label{EqCocyForTransis}
g_{kj}(u)\ g_{ji}(u)=g_{ki}(u),\qquad \zeta_{kj}(u,g_{ji}(u))\ \zeta_{ji}(u,g)=\zeta_{ji}(u,g)
\end{eqnarray}
on threefold intersections $U_{ijk}\ni u$.
Conversely, given families $\{g_{ji}\}, \{\zeta_{ji}\}$ 
which satisfy (\ref{EqCocyForTransis}) we (re-) obtain
a pair by the usual quotients.
The two families $\{g_{ji}\},\{\zeta_{ji}\}$ of transition functions
can also be considered as one single family
\begin{eqnarray*}
g_{ji}\times\zeta_{ji}: U_{ji}\to G\ltimes C(G,\U(1)),
\end{eqnarray*}
where $G\ltimes C(G,\U(1))$ is the semi-direct product.
Note that the (topological) group
$${\rm A_{\rm Par}}:= G\ltimes C(G,\U(1))$$
is the automorphism group of the trivial pair over the point
$G\times \U(1)\to G\to\ast$,
where $(g,f)\in {\rm A_{\rm Par}}$ acts from the left 
on $(h,z)\in G\times \U(1)$ by 
$(g,f)\cdot(h,z):=(g+h, f(h)z)$.

Let $P\to B$ be a ${\rm A_{\rm Par}}$-principal fibre bundle.
We obtain a pair over $B$ by associating the trivial pair 
over the point\footnote{
For a group $H$ acting on a space $X$ from the right 
and on a space $Y$ from the left 
we denote by $X\times_{H}Y$ the quotient of $X\times Y$ by the induced
right action $(x,y)\cdot h:= (x\cdot h, h^{-1}\cdot y)$ as usual.
}
\begin{eqnarray*}
	F_P:= P\times_{{\rm A_{\rm Par}}}( G\times \U(1)),\quad 
	E_P:= P\times_{{\rm A_{\rm Par}}} G.
\end{eqnarray*}
A morphism of ${\rm A_{\rm Par}}$-bundles over $B$
\begin{eqnarray*}
	\xymatrix{
	P\ar[r]\ar[d]&P'\ar[d]\\
	B\ar[r]^=&B
	}
\end{eqnarray*}
induces a morphism of pairs over $B$
\begin{eqnarray*}
	\xymatrix{
	F_P\ar[r]\ar[d]&F_{P'}\ar[d]\\
	E_P\ar[r]\ar[d]&E_{P'}\ar[d]\\
	B\ar[r]^=&B
	},
\end{eqnarray*}
i.e. we have constructed a functor 
 from the category of 
 ${\rm A_{\rm Par}}$-principal fibre bundles over $B$ 
 to the category of pairs over $B$.

\begin{prop}\label{PropEqOfCatsPairsAndBunAFull}
	The functor
	\begin{eqnarray*}
		(P\to B)\mapsto (F_P\to E_P\to B)
	\end{eqnarray*}
	 is an equivalence of categories.
\end{prop}

\begin{pf}
	We construct a functor in opposite direction.
	Let $F\stackrel{q}{\to} E\to B$ be a pair. 
	For $b\in B$ we define $P_b$ to be the set of 
	tupels $(e,s)$ of arbitrary elements $e\in E|_b$ and
	subsets $s\subset F|_{E|_b}$ with the property
	that $  q|_{s}: s\to E|_b$ is a homoeomorphism.
	Denote by $s(x)$ the unique element in $s$
	such that $q(s(x))=x$, for $x\in E|_b$, so 
	$s=\{s(x)| x\in E|_b\}$.
	The set $P_b$ is a ${\rm A_{\rm Par}}$-torsor. 
	In fact, for $(g,f)\in {\rm A_{\rm Par}}$ we have a
	right action on $(e,s)\in P_b$ by
	\begin{eqnarray}\label{EqNetteActionImZugNachWien}
	(e,s)\diamond\hspace{-0.242cm}\cdot\ (g,f):= (e\cdot g, s\diamond_e(g,f)),
	\end{eqnarray}
	where  $e\cdot g$ is the principal action,
	and $s\diamond_e(g,f)$ is the subset 
	$\{ s(e\cdot (h+g))\cdot f(h)| h\in G \}$.
	As (the graph of) a function $f:G\to\U(1)$ 
	is a subset of $G\times\U(1)$, an atlas $U_i$ 
	of the pair induces ${\rm A_{\rm Par}}$-equivariant bijections 
	$\varphi_i: \coprod_{b\in U_i} P_b\cong U_i\times G\ltimes C(G,\U(1))$.
	We equip 
	$\coprod_{b\in B} P_b$ with the coarsest topology 
	such that all $\varphi_i$ are homoeomorphisms.
	This topology is independent of the chosen atlas.
	The canonical map $\coprod_{b\in B} P_b\to B$ is now
	an ${\rm A_{\rm Par}}$-principal fibre bundle.
	It is rather obvious that a morphism of pairs
	induces a morphism of ${\rm A_{\rm Par}}$ principal bundles,
	and it is straight forward to verify that the compositions 
	of the two functors in game are naturally isomorphic to 
	the identity functors.  
\end{pf}

We give two simple but important examples of pairs
which correspond to ${\rm A_{\rm Par}}$-bundles which have 
reductions of the structure group ${\rm A_{\rm Par}}$ to two 
distinguished  subgroups.

First, bundles which have a reduction to the subgroup 
$G\subset {\rm A_{\rm Par}}$ lead to 
pairs isomorphic to
$$
\xymatrix{
E\times \U(1)\ar[d]\\
E\ar[d]\\ 
B
}
$$
for a $G$-bundle $E\to B$.

Second, we have $\widehat G \subset C(G,\U(1))\subset {\rm A_{\rm Par}}$.
Bundles which admit reductions to the subgroup
$\widehat G\subset {\rm A_{\rm Par}}$ lead to pairs 
 of the form
\begin{equation}\label{diagTheGenRingPair}
\xymatrix{
F_{\widehat E}\ar[d]\\
B\times G\ar[d]\\ 
B
}
\end{equation} 
where
$
F_{\widehat E}:= \widehat E\times_{\widehat G} (G\times \U(1)),
$
for a $\widehat G$-bundle $\widehat E\to B$.

The latter example of a pair and its 
algebraic properties are the contend 
of the next section.

\section{Ring pairs}
\label{SecRingPairs}

The trivial pair over the point admits 
a multiplication in the sense that 
\begin{eqnarray}\label{DiagTheTrivialRingPairInParis}
\xymatrix{
(G\times \U(1))\times (G\times\U(1))
\ar[rr]^{\qquad+\ \times\ \cdot}
\ar[d]&
&G\times\U(1)\ar[d]\\
G\times G\ar[rr]^+\ar[d]&&G\ar[d]\\ 
\ast\ar[rr]^=&&\ast
}
\end{eqnarray}
commutes, where $+\times\cdot$ is the map
${((g,t),(h,s))\mapsto (g+h,ts)}$.
We take this as the local model for the notion 
of ring pairs which we introduce next.
Consider an arbitrary pair $F\to E\to B$ and 
the following diagram of pullbacks.
$$
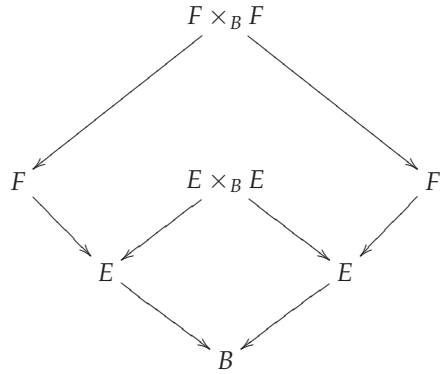
\resizebox*{!}{5cm}{
\xymatrix{
&&F\times_B  F\ar[rdrd]\ar[ldld]&&\\
&&&&\\
F\ar[rd]&&E\times_B  E\ar[rd]\ar[ld]&&F\ar[ld]\\
&E\ar[rd]&&E\ar[ld]&\\
&&B&&
}}
$$
The canonical map
$F\times_B  F\to E\times_B  E$
is a $\U(1)\times \U(1)$-bundle and
$E\times_B  E\to B$
is a $G\times G$-bundle.

\begin{defi}
	A {\bf ring pair} is a pair
	$F\to E\to B$ together with 
	a commutative diagram
	\begin{equation}\label{diagDefiofRingPair}
		\xymatrix{
		F\times_B F\ar[r]\ar[d]& F\ar[d]\\
		E\times_B E\ar[r]\ar[d]&E\ar[d]\\ 
		B\ar[r]^=&B
		}
	\end{equation}
	such that there exists an atlas with the property that for each
	chart $U_i$ and  $b\in U_i$ the restriction 
	$$
	\xymatrix{
	F|_{E|_b}\ar[r]\ar[d]& G\times \U(1)\ar[d]\\
	E|_b\ar[r]\ar[d]& G\ar[d]\\
	b\ar[r]& b
	}
	$$
	of the chart induces 
	$$
	\xymatrix{
	(G\times \U(1))\times(G\times\U(1))\ar[rr]^{\qquad+\ \times\ \cdot}\ar[d]&
	&G\times\U(1)\ar	[d]\\
	G\times G\ar[rr]^+\ar[d]&&G\ar[d]\\ 
	b\ar[rr]^=&&b
	}
	$$
	from the restriction (\ref{diagDefiofRingPair})$|_b$.

	The horizontal maps in (\ref{diagDefiofRingPair})
	are called  the {\bf multiplication maps} of the pair, 
	and an atlas satisfying the conditions above is
	called a {\bf ring atlas} for the ring pair.
\end{defi}

An example of a ring pair can be made out of the pair 
$F_{\widehat E}\to B\times G\to B$ from (\ref{diagTheGenRingPair}).
Namely, this pair admits canonical maps
\begin{eqnarray}\label{EqGetMyFirstRingStr}
\resizebox{!}{0,55cm}{$
\begin{array}{rclrcl}
F_{\widehat E}\times_B F_{\widehat E}&\to&F_{\widehat E}
&(B\times G)\times_B (B\times G)&\to &B\times G\\
([\hat e,g,t],[\hat e\cdot\chi,h,s])&\mapsto&[\hat e,g+h,\langle h,\chi\rangle t s ]&
((b,g),(b,h))&\mapsto&(b,g+h)
\end{array}
$}
\end{eqnarray}
which give this pair the structure of a ring pair. 
A ring atlas is obtained from
the local trivialisations of the bundle $\widehat E\to B$.
We will see below, that all examples are of this form.

\begin{rem}\label{RemRingPairIsGroupBndl}
For a ring pair $F\to E\to B$ the map $F\to B$ 
is a bundle of groups with fibre isomorphic to $G\times \U(1)$.

In particular,  the multiplication map $\mu$ is associative and commutative,
and there is a canonical section $\sigma: B\to F$.
If $\mu:F\times_B F\to F$ denotes the multiplication map
and $x\in F$ is an element over $b\in B$, then there
is a unique $x^\dag\in F$ over $b$ with the property
\begin{eqnarray*}\label{EqTheInversOfMulti}
\mu(x, x^\dag)=\mu(x^\dag,x)=\sigma(b).
\end{eqnarray*}
$^\dag:F\to F$ is anti-linear, i.e.
$
(x\cdot t)^\dag=x^\dag\cdot\overline t,
$ 
for the complex conjugate $\overline t$ of $t\in\U(1)$.
\end{rem}

\begin{pf}
The associativty and commutativity of $\mu$ follows 
as the multiplication map is locally associative and commutative.

To define the section $\sigma:B\to F$ we define 
$\sigma(b)\in F|_{E|_b}$ to be the pre-image of $(0,1)$
under any chart $F|_{E|_b}\to G\times \U(1)$ of an ring atlas.
We have to check that this is well-defined.
So let $\{U_i\}$ be a ring atlas with transition 
functions $g_{ij}, \zeta_{ij}$ satisfying (\ref{EqCocyForTransis}).
Assume $b\in U_i$ and $b\in U_j$, then we have a commutative diagram 
\begin{eqnarray*}\label{DiagChangeOfChart}
\xymatrix{
(U_i\times G\times\U(1))\times_{U_i}( U_i\times G\times\U(1))&&\\
\ar@{^(->}[u]^b \ar[d]_{\underset{(g_{ji}(b)+g,g_{ji}(b)+h,\zeta_{ji}(b,g)t, \zeta_{ji}(b,h) s)} 
{\overset{(g,h,t,s)}{\rotatebox{270}{$\mapsto$}}}}
 G\times G \times \U(1)\times \U(1)
\ar[rr]^{\qquad\ +\ \times\ \cdot }&&
G\times \U(1)\ar[d]^{\underset{(g_{ji}(b)+g,\zeta_{ji}(b,g)t)} 
{\overset{(g,t)}{\rotatebox{270}{$\mapsto$}}}}\\
\ar@{_(->}[d]_b  G\times G\times \U(1)\times\U(1)\ar[rr]^{\qquad\ +\ \times\ \cdot}&&
G\times \U(1)\\
(U_j\times G\times\U(1))\times_{U_j}( U_j\times G\times\U(1))&&
}.
\end{eqnarray*}
The commutativity of this diagram implies that
\begin{eqnarray}\label{EqTrasiForRingPairss}
g_{ji}(b)=0 \text{ and }\zeta_{ji}(b,g)\zeta_{ji}(b,h)=\zeta_{ji}(b,g+h),
\end{eqnarray}
so in particular $\zeta_{ji}(b,0)=1$.
Thus, the element $(0,1)\in G\times \U(1)$ is fixed by 
any change of charts. It follows that $\sigma(b)$ is well-defined.

Fix an $x\in F$ over $b\in B$. The map 
$\mu(x,\_): F|_{E|_b}\to F|_{E|_b}$ is an isomorphism 
thus there exists a unique $ x^\dag\in F$ such that
$\mu(x, x^\dag)=\mu( x^\dag,x)=\sigma(b).$
The anti-linearity of $^\dag$ is clear then.
\end{pf}
 
The composition $B\stackrel{\sigma}{\to}F\to E$
is a section of $E\to B$, hence there is a 
canonical trivialisation
\begin{eqnarray*}
E\stackrel{\cong}{\to} B\times G.
\end{eqnarray*}
By this trivialisation,  $E\to B$ has a canonical 
structure of a bundle of groups.
A bundle of groups consists of the same data as a 
groupoid where source and target maps are equal.
Therefore the $\U(1)$-bundle $F\to E$ is
a $\U(1)$-central extension of groupoids
\begin{eqnarray}\label{DiagCentrExtOfGrpoid}
\xymatrix{
B\!\times\!\U(1)\ \ar@<2pt>[d]\ar@<-2pt>[d]\ar@{>->}[r]&F\ar@{->>}[r]\ar@<2pt>[d]\ar@<-2pt>[d] 
&E\ar@<2pt>[d]\ar@<-2pt>[d]\\
B\ar[r]^=&B\ar[r]^=&B
}.
\end{eqnarray}

\begin{rem}\label{RemOnAlgOfRP}
	In section \ref{SecTheFourierTransform}
	we will be concerned about the $C^*$-algebra
	$C^*(E,F)$ of a ring pair $F\to E\to B$. 	
	By thinking of ring pairs as central extensions
	it is possible to give a quick definition of $C^*(E,F)$. 
	It is just the $C^*$-algebra of the central 
	extension (\ref{DiagCentrExtOfGrpoid}).
\end{rem}

Equality (\ref{EqTrasiForRingPairss}) determines the 
possible transition functions in a ring atlas.
There exist  $\chi_{ji}:U_{ji}\to \widehat G$ 
such that $\zeta_{ji}(u,g)=\langle g,\chi_{ji}(u)\rangle$.
Thus, in a ring atlas the transition functions
are maps
\begin{eqnarray*}
U_{ji}\to \widehat G\subset {\rm A_{\rm Par}}.
\end{eqnarray*}
Note that $\widehat G$ is the 
automorphism group of the trivial ring pair 
over the point (\ref{DiagTheTrivialRingPairInParis})
in the category of ring pairs.
Here a morphism in the category of ring pairs is  a morphisms of pairs 
such that the induced diagram of multiplication maps 
commutes.
The category of ring pairs and the
category of $\widehat G$-bundles 
are proper  subcategories (i.e. \emph{not} full subcategories) of
the categories of ring pairs and 
${\rm A_{\rm Par}}$-bundles.
Therefore the following proposition is not just
a corollary of Proposition 
\ref{PropEqOfCatsPairsAndBunAFull}.
However, the bundle theoretic argument  
used in the proof is quite the same.

\begin{prop}\label{PropMyFirstEqOfCats}
	The functor
	\begin{eqnarray}\label{EqEqOfCatstevens}
		\big(\widehat E\to B\big)\mapsto 
		\big(F_{\widehat E}\to B\times G\to B,
		(\ref{EqGetMyFirstRingStr})\big)
	\end{eqnarray}
	is an equivalence of categories from the category 
	of $\widehat G$-bundles over $B$
	to the category of ring pairs over $B$.
\end{prop}

\begin{pf}
We define a $\widehat G$-bundle from the data of a ring pair
$F\to E\to B$. Let $\sigma, \mu$ be as in 
Remark \ref{RemRingPairIsGroupBndl}, and
let $r:F\to B\times G$ denote the composition  $F\to E\cong B\times G$.
For $b\in B$, let $\widehat E_b$ be the set of all subsets $\hat e\subset F|_{E|_b}$
with the property that, first, $r|_{\hat e}:\hat e\to b\times G$  is a homoeomorphism,
and second, $ \mu( \hat e(g),\hat e(h))=\hat e(g+h)$. Here we used the notation
$\hat e(g):=(r|_{\hat e})^{-1}(b\times g)\in \hat e$.
The set $\widehat E_b$ is a $\widehat G$-torsor subject to the action
$$
\hat e\diamond \chi:=\{ \hat e(g)\cdot\langle g,\chi\rangle\in F| g\in G\},
$$
where $\U(1)$ acts on $F$ by the given principal action.
Let $\widehat E:=\bigsqcup_{b\in B} \widehat E_b$.
We give $\widehat E$ a topology such that
the canonical map $\widehat E\to B$ becomes a
$\widehat G$-principal bundle.
Namely, if $U_i$ is any ring atlas for the pair
$F\to E\to B$, then 
$U_i\times G\times \U(1)\cong F|_{E|_{U_i}}$
identifies the characters $\chi\in \widehat G$ 
understood as graphs $\chi \subset G\times\U(1)$ 
with the elements of 
$\bigsqcup_{b\in U_i} \widehat E_b\subset\widehat E$. 
We take the coarsest topology on
$\widehat E$ such that all the bijections
$\bigsqcup_{b\in U_i} \widehat E_b\cong U_i\times \widehat G$ 
are homoeomorphisms.

A morphism of ring pairs preserves the multiplication maps,
it induces a morphism between the constructed 
$\widehat G$-bundles. 
So we have defined a functor 
from ring pairs to $\widehat G$-bundles.

To check that the compositions of the two functors 
are naturally isomorphic to the corresponding 
identity functors is a tedious manipulation.
\end{pf}

We give an example of a (non-trivial) ring pair $F\to E\to B$
where both of the bundles $F\to E$
and $E\to B$ are trivial, but the corresponding 
$\widehat G$-bundle $\widehat E$ is not,
i.e. the non-triviality is hidden in the multiplication maps.

\begin{exa}\label{ExaNTRPWUTP}
Take $G:={\rm O}(1)\cong \widehat G$ the group with two elements,
and let $B:=\R {\rm P}^1$ be the 1-diminsional projective space.
The quotient map $S^1\to \R {\rm P}^1$ is a 
non-trivial $\widehat{G}$-bundle, but the 
$\U(1)$-bundle $F_{S^1}\to \R {\rm P}^1\times {\rm O}(1)$
is trivialisable in the category of 
$\U(1)$-bundles over $\R {\rm P}^1\times {\rm O}(1)$
since $H^2(\R {\rm P}^1\times{\rm O}(1))=0$.
Thus, the choice of a trivialisation
$\R {\rm P}^1\times{\rm O}(1)\times\U(1)\cong F_{S^1}$
and the multiplication map on $F_{S^1}$ yield
a non-trivial ring pair with an underlying trivial pair.
\end{exa}

\section{Module pairs}
\label{SecModulePairs}
\noindent
Having the notion of ring pairs at hand we can define the 
notion of a module pair over a ring pair. 

\begin{defi}\label{DefinModulePair}
A {\bf module pair} over a ring pair
$F_0\to E_0\to B$ is a pair $F\to E\to B$
together with 
a commutative diagram
\begin{equation}\label{diagDefiofModPair}
\xymatrix{
F\times_B F_0\ar[r]\ar[d]& F\ar[d]\\
E\times_B E_0\ar[r]\ar[d]&E\ar[d]\\ 
B\ar[r]^=&B
}
\end{equation}
such that there exists an atlas for $F\to E\to B$
and a ring atlas for $F_0\to E_0\to B$ over the
same open covering $\{U_i\}$ of $B$
with the property that for each chart of the ring 
atlas over $U_i$ there is a chart  of $F\to E\to B$
over $U_i$ such that for all $b\in U_i$ 
the restrictions  
$$
\xymatrix{
F|_{E|_b}\ar[r]\ar[d]& G\times \U(1)\ar[d]&
F_0|_{E_0|_b}\ar[r]\ar[d]& G\times \U(1)\ar[d]
\\
E|_b\ar[r]\ar[d]& G\ar[d]&
E_0|_b\ar[r]\ar[d]& G\ar[d]&
\\
b\ar[r]&b&
b\ar[r]&b
}
$$
of these charts induce 
$$
\xymatrix{
(G\times \U(1))\times (G\times\U(1))\ar[rr]^{\qquad+\ \times\ \cdot}\ar[d]&&G\times\U(1)\ar[d]\\
G\times G\ar[rr]^+\ar[d]&&G\ar[d]\\ 
b\ar[rr]^=&&b
}
$$
from the restriction (\ref{diagDefiofModPair})$|_b$.

An  atlas as above is called a {\bf module atlas}
for $F\to E\to B$. 
The horizontal maps in (\ref{diagDefiofModPair})
are called  {\bf action maps}.
\end{defi}

\begin{rem}\label{RemAboutModPairs}
\begin{enumerate}
\item
The multiplication maps of a ring pair and the action
maps of a module pair are associative. I.e. 
if $\mu:F_0\times_B F_0\to F_0$ and
$\varrho:F\times_B F_0\to F$ 
denote the multiplication map and the action 
map  of a ring pair and a module pair,
then
\begin{eqnarray*}
\varrho(x,\mu(y_0,z_0))=\varrho(\varrho(x,y_0),z_0)
\end{eqnarray*}
holds.
\item
Fix  $x\in F$ in the fibre over $b\in B$. The map
$\varrho(x,\_): F_0|_{E_0|_b}\to F|_{E|_b}$ is 
an isomorphism, thus 
there is a well-defined map 
$\sigma_\varrho: F\times_B F\to F_0$
unique with respect to the property
\begin{eqnarray*}
\varrho(x,\sigma_\varrho(x,y))=y.
\end{eqnarray*}
One should think of $\sigma_\varrho$ as
an $ F_0$-valued ``scalar product''
as it is linear and $F_0$-linear in the second entry, i.e.
\begin{eqnarray*}
\sigma_\varrho(x,y\cdot t)=\sigma_\varrho(x,y)\cdot t,\quad \sigma_\varrho(x,\varrho(y,z_0)=\mu(\sigma_\varrho(x,y), z_0),
\end{eqnarray*}
it is anti-symmetric, i.e.
\begin{eqnarray*}
\sigma_\varrho(y,x)={\sigma_\varrho(x,y)^\dag},
\end{eqnarray*}
and it is positive definite, i.e.
\begin{eqnarray*}
\sigma_\varrho(x,x)=\sigma(b).
\end{eqnarray*}
Here $\sigma: B\to F_0$  and ${{\ \ }^\dag}:F_0\to F_0$ 
are as in Remark \ref{RemRingPairIsGroupBndl}.
\end{enumerate}
\end{rem}

Let us return to the groupoid point of view. 
As explained in diagram (\ref{DiagCentrExtOfGrpoid}), 
we can think of a ring pair $F_0\to E_0\to B$ as of 
a central extension of groupoids
\begin{eqnarray*}
\xymatrix{
B\!\times\!\U(1)\ \ar@<2pt>[d]\ar@<-2pt>[d]\ar@{>->}[r]&F_0\ar@{->>}[r]\ar@<2pt>[d]\ar@<-2pt>[d] 
&E_0\ar@<2pt>[d]\ar@<-2pt>[d]\\
B\ar[r]^=&B\ar[r]^=&B
}.
\end{eqnarray*}
Let $E\to B$ be a $G$-principal bundle.
It can be thought of as a Morita self-equivalence 
\begin{eqnarray*}
	\xymatrix{
	E_0 &
	\ar@{}[l]|{\ \ \resizebox{!}{0,6cm}{\rotatebox {90}{{\rotatebox {90}{
	{\rotatebox {90}{$\circlearrowright$}}}}}}}
	E\ar[dl]\ar[dr] \ar@{}[r]|{\resizebox{!}{0,6cm}{\rotatebox {90}{$\circlearrowleft$}}}
	& E_0\\
	B\ar@{<-}@<2pt>[u]\ar@{<-}@<-2pt>[u]&&B\ar@{<-}@<2pt>[u]\ar@{<-}@<-2pt>[u]
	}.
\end{eqnarray*}
where the left and the right action of $E_0\rightrightarrows B$ coincide
and commute by the abelianess of $G$.
By the Remark \ref{RemAboutModPairs}, 
a module pair  $F\to E\to B$ over
$F_0\to E_0\to B$ gives rise to a $\U(1)$-equivariant  
Morita self-equivalence
\begin{eqnarray}\label{DiaU1EqME}
	\Big(
	\xymatrix{
	F_0 &
	\ar@{}[l]|{\ \ \resizebox{!}{0,6cm}{\rotatebox {90}{{\rotatebox {90}{
	{\rotatebox {90}{$\circlearrowright$}}}}}}}
	F\ar[dl]\ar[dr] \ar@{}[r]|{\resizebox{!}{0,6cm}{\rotatebox {90}{$\circlearrowleft$}}}
	& F_0\\
	B\ar@{<-}@<2pt>[u]\ar@{<-}@<-2pt>[u]&&B\ar@{<-}@<2pt>[u]\ar@{<-}@<-2pt>[u]
	}\Big)\ {\rotatebox {90}{$\circlearrowleft$}}\ \U(1).
\end{eqnarray}

\begin{rem}
	Viewing to module pairs as Morita equivalences 
	of groupoids will be helpful in section \ref{SecTheFourierTransform}.
	There we will define the Hilbert $C^*$-module $H(E,F)$ of a module pair
	$F\to E\to B$. The shortest way  to  do this is to define $H(E,F)$ as 
	the Hilbert $C^*$-module 
	of the  $\U(1)$-equivariant  Morita equivalence (\ref{DiaU1EqME}). 
	It is a module over the  $C^*$-algebra $C^*(E_0,F_0)$ 
	of the ring pair $F_0\to E_0\to B$
	(s. Remark \ref{RemOnAlgOfRP}).
\end{rem}

Let $F\to E\to B$ be a pair
and let $\widehat E\to B$ be a $\widehat G$-bundle.
It is natural to ask whether or not one can turn
$F\to E\to B$ into a module pair over 
$F_{\widehat E}\to B\times G\to B$.
To answer this question consider first the cup product
$$
\cup : \check H^1(B,\underline G)\times \check H^1(B,\underline{\widehat G}) 
\to\check H^2(B,\underline\U(1)),
$$
i.e. if $g_{ji}:U_{ij}\to G$ and $\chi_{ji}:U_{ij}\to\widehat G$  
are transition functions of the bundles  $E$ and $\widehat E$, then
the class $[E]\cup [\widehat E]$ is represented by the cocycle 
$\alpha_{kji}:U_{ijk}\to \U(1)$ with 
\begin{eqnarray}\label{EqTheCocDefForCupPro}
\alpha_{kji}(u):=\langle g_{ij}(u),\chi_{jk}(u)\rangle\in\U(1).
\end{eqnarray}
Assume that the class $[E]\cup [\widehat E]$ vanishes, 
so  let $s_{ji}:U_{ji}\to \U(1)$ be such that 
$\check\delta\{s_{ji}\}=\{\alpha_{kji}\}$. Then we can define
\begin{equation}\label{EqDasGibtNenTorsorAlter}
\zeta_{ji}(u,g):=\langle g,\chi_{ji}\rangle s_{ji}(u)^{-1},
\end{equation}
and it is immediate that $g_{ji}$ and $\zeta_{ji}$ satisfy 
(\ref{EqCocyForTransis}), i.e. $\zeta_{ij}$ is a family of
transition functions for a $\U(1)$-bundle over $E$.
As $s_{ji}$ is unique  only up to a 1-cocycle in $\check Z^1(\{U_i\},\underline{\U(1)})$
we see that the class of $\zeta_{ji}$ in $\check H^1(E,\underline{\U(1)})$
is only determined up to the natural action 
of $\check H^1(B,\underline{\U(1)})$ on
$\check H^1(E,\underline{\U(1)})$ given by pullback
$p^*:\check H^1(B,\underline{\U(1)})\to \check H^1(E,\underline{\U(1)})$
along the projection $p:E\to B$.
In other words, if $[E]\cup[\widehat E]=0$, 
then (\ref{EqDasGibtNenTorsorAlter}) 
defines a subset
\begin{eqnarray}\label{EqAllesSenkrecht}
[\widehat E]^\perp\subset \check H^1(E,\underline{ \U(1)})
\end{eqnarray}
which is a $p^*(\check H^1(B,\underline{\U(1)})$-torsor.

\begin{prop}
	$F\to E\to B$ is a module pair over 
	$F_{\widehat E}\to B\times G\to B$ if and only if 
	$[E]\cup [\widehat E]=0$ and 
	$[F]\in [\widehat E]^\perp$.
\end{prop}

\begin{pf}
	Let $F\to E\to B$ be a module pair. 
	Let $U_i$ be a module atlas with transition
	functions $g_{ji},\zeta_{ji}$. 
	Let $\zeta_{ji}^0$ denote the transition functions 
	in the ring atlas, by equality (\ref{EqTrasiForRingPairss}) 
	we have $\zeta^0_{ji}(b,h)=\langle h,\chi_{ji}(b)\rangle$.
	Similar to equality (\ref{EqTrasiForRingPairss}) we
	obtain that 
	$\zeta_{ji}(b,g) \zeta^0_{ji}(b,h) =\zeta_{ji}(b,g+h)$.
	It follows that $\zeta_{ji}(b,h)=\langle h,\chi_{ji}(b)\rangle s_{ji}(b)^{-1}$,
	for $s_{ji}(b):=\zeta_{ji}(b,0)^{-1}$. One computes
	$s_{kj}(b)s_{ki}(b)^{-1}s_{ji}(b)=\langle g_{ji}(b),\chi_{kj}(b)\rangle$.
	Hence $[E]\cup[\widehat E]=0$ and $[F]\in[\widehat E]^\perp.$
	
	Conversly, let $F\to E\to B$ be a pair such that 
	$[E]\cup[\widehat E]=0$ and $[F]\in[\widehat E]^\perp.$
	Take a ring atlas for the ring pair and an atlas of the pair 
	such that the transition functions 
	satisfy (\ref{EqDasGibtNenTorsorAlter}). Then we define 
	the action maps by demanding that the chosen atlas
	is a module atlas. It is straight forward to check that this
	is well-defined.	
\end{pf}

Fix two bundles $E\to B$ and $\widehat E\to B$ with
$[E]\cup [\widehat E]=0$.
The set $[\widehat E]^\perp$ is the set of isomorphism 
classes of $\U(1)$-bundles over $E$ which admit
a module structure, but the isomorphisms do not
have to be compatible with the module structure.

Let ${\rm Mod}(E,\widehat E)$ be the groupoid 
of module pairs $F\to E\to B$ over 
$F_{\widehat E}\to B\times G\to B$,
where a morphism from $F\to E\to B$ to $F'\to E\to B$
is a bundle morphism over $E$
\begin{eqnarray}
\xymatrix{
F\ar[r]^\cong\ar[d]& F'\ar[d]\\
E\ar[r]^=& E
}
\end{eqnarray}
such that the induced diagram 
\begin{eqnarray}\label{DiagIsoModPairs2}
\xymatrix{
F\times_B F_{\widehat E}\ar[d]\ar[r]& F\ar[d]\\
F'\times_B F_{\widehat E}\ar[r]& F'
}
\end{eqnarray}
commutes. Here the horizontal maps are the action maps of the module pairs.

Let $[{\rm Mod}(E,\widehat E)]$ be the set of isomorphism classes 
of ${\rm Mod}(E,\widehat E)$.
By pullback, there is a well defined, natural action of
$\check H^1(B,\underline {\U(1)})$ on $[{\rm Mod}(E,\widehat E)]$. 
In fact, let $\eta$ and $\eta'$ be two $\U(1)$-bundles which represent the same class 
$[\eta]=[\eta']\in \check H^1(B,\underline{\U(1)})$,  and let
$F\to E\stackrel{p}{\to} B$ be a module pair. Then
$F\otimes p^*\eta$ and $F\otimes p^*\eta'$ become module pairs in the obvious way,
and they are isomorphic in ${\rm Mod}(E,\widehat E)$.

\begin{prop}\label{PropIsoClassesOfMod}
The set $[{\rm Mod}(E,\widehat E)]$ is a $\check H^1(B,\underline {\U(1)})$-torsor.
\end{prop}

\begin{pf}
The action is transitive: 
If $[F]$ and $[F']$ are classes in $[{\rm Mod}(E,\widehat E)]$, 
then it follows from (\ref{EqDasGibtNenTorsorAlter}) that there exists a $\U(1)$-bundle $\eta\to B$
such that $F$ and $F'\otimes p^*\eta$ are isomorphic in ${\rm Mod}(E,\widehat E)$.

The action is free:
Let $\{U_i\}$ be a sufficiently refined open cover of $B$ such that
the following algebraic arguments are meaningful.
Let $\eta\to B$ be a bundle with transition functions $\eta_{ji}:U_i\to\U(1)$ 
such that $F$ and $F\otimes p^*\eta$ 
represent the same class in ${\rm Mod}(E,\widehat E)$.
In particular, as they are isomorphic over $E$ and as their
transition functions vary by $\eta_{ji}$,  there
are local isomorphisms $\varphi_i:U_i\times G\to \U(1)$ such that
$\varphi_j(u,g_{ji}(u)+g)^{-1}\ \varphi_i(u,g)=\eta_{ji}(u)$.	
Using the commutativity of (\ref{DiagIsoModPairs2}),
it is easy to see that each $\varphi_i$
is constant in the $G$-argument, so
$\varphi_j(u,0)^{-1}\ \varphi_i(u,0)=\eta_{ji}(u)$, 
i.e. the class of $\eta_{ji}$ is trivial in $\check H^1(B,\U(1))$.
\end{pf}

Now, fix two bundles $E\to B$ and $\widehat E\to B$ without any further
assumption about  the class $[E]\cup [\widehat E]\in \check H^2(B,\underline{\U(1)})$.
For any open set $U\subset B$ we can consider the groupoid 
${\rm Mod}_{E,\widehat E}(U):= {\rm Mod}(E|_U,\widehat E|_U)$, 
where $E|_U\to U$ and $\widehat E|_U\to U$ are the
pulbacks of $E$ and $\widehat E$ along the inclusion $U\subset B$.

\begin{prop}\label{PropDieModGerbe}
The assignment 
\begin{eqnarray}\label{EqDieModGerbe}
{\rm Mod}_{E,\widehat E} : U\mapsto {\rm Mod}(E|_U,\widehat E|_U)
\end{eqnarray}
is a $\U(1)$-banded gerbe over $B$, and its class given by
$[E]\cup [\widehat E]\in \check H^2(B,\underline{\U(1)})$.
\end{prop}

\begin{pf}
Is is immediate that (\ref{EqDieModGerbe}) is a stack.
It is also clear from the definition of module pairs 
that (\ref{EqDieModGerbe}) is locally non empty, and that
locally any two objects are isomorphic. 
Hence (\ref{EqDieModGerbe}) is a gerbe.

Let $F\in {\rm Mod}(E|_U,\widehat E|_U)$ be an object,
and let $V\subset U$ be open.
Pullback along $E|_U|_V\to V$ defines an injection
$C(V,\U(1))\hookrightarrow C(E|_U|_V,\U(1))\hookrightarrow {\rm Aut}(F|_{E|_U|_V})$. 
The commutativity of (\ref{DiagIsoModPairs2}) implies 
that every automorphism of $F$ is of this form. 
Thus for each open $U$ the sheaf of automorphisms 
$\underline{\rm Aut}(F)$ on $U$ is isomorphic to $C(\_,\U(1))$.
Hence the gerbe (\ref{EqDieModGerbe}) is $\U(1)$-banded.

Let $U_i$ be an open cover together with trivialisations  
$E|_{U_i}\to U_i\times G$ and $\widehat E|_{U_i}\to U_i\times\widehat G$
and transition functions $g_{ij}:U_{ij}\to G$ and $\chi_{ij}:U_{ij}\to \widehat G$.
We choose for each $i$ an $F_i\in {\rm Mod}(E|_{U_i},\widehat E|_{U_i})$,
namely $F_i:= E|_{U_i}\times \U(1)$. 
The trivialisations of $E|_{U_i}$ and $\widehat E|_{U_i}$
define a unique module structure on $F_i\to E|_{U_i}\to U_i$
over $F_{\widehat E|_{U_i}}\to U_i\times G\to U_i$.
Let $F_i|_ j\in {\rm Mod}(E|_{U_{ij}},\widehat E|_{U_{ij}})$
denote the pullback of $F_i$ along the composition
$E|_{U_{ij}}\cong E|_{U_{i}}|_{U_{ij}}\hookrightarrow 
 E|_{U_{i}}$.
 Interchanging $i$ and $j$ we  also find
 $F_j|_i\in {\rm Mod}(E|_{U_{ij}},\widehat E|_{U_{ij}})$, 
and we define a morphism $f_{ji}$ between these by 
$$
\xymatrix{
F_{i}|_{j}\ar[d]\ar[rrrr]^{f_{ji}}&&&&F_j|_i\\
U_{ji}\times G\times\U(1)
\ar[rrrr]_{ (u,g,z)\mapsto (u,g+g_{ji}(u),\langle g,\chi_{ji}(u)\rangle z)   }
&&&&U_{ji}\times G\times \U(1)\ar[u]
},
$$
where the vertical trivialisations are given by the 
corresponding trivialisations over $U_i, U_j$ respectively,
Then on threefold intersections we have the automorphism
$f_{ki}^{-1}\circ f_{kj}\circ f_{ji}\in\underline{\rm Aut}(F_i)(U_{ijk})$ which 
is mapped to a $\check C$ech cocycle $ \beta_{kji}\in C(U_{ijk},\U(1)) $ 
under the identification $\underline{\rm Aut}(F_i)(U_{ijk})\cong C(U_{ijk},\U(1))$.
This cocycle represents the class of the gerbe. We have
$\beta_{kji}(u)= \langle g_{ji}(u),\chi_{kj}(u)\rangle$, and 
if we compare with the cocycle  (\ref{EqTheCocDefForCupPro})
we find $\alpha_{kji}=\beta_{kji}$.
Hence the class of  the gerbe (\ref{EqDieModGerbe}) is $[E]\cup[\widehat E]$.
\end{pf}

\section{Pontrjagin duality triples}
\noindent
To introduce the main objects of our interest 
we first define their local model.

\begin{defi}\label{defitDualityTripleOverPt}
The {\bf trivial Pontrjagin duality triple}  over the point  is 
the diagram 
\begin{equation}\label{diagDefiTrivTdual}
\resizebox*{!}{4cm}{
\xymatrix{
&G\times\widehat G\times \U(1)\ar[rd]\ar[ld]\ar[rr]^{\pi}&
&G\times\widehat G \times\U(1)\ar[rd]\ar[ld]&\\
G\times\U(1)\ar[rd]&&G\times \widehat G\ar[rd]\ar[ld]&
&\widehat G\times\U(1)\ar[ld] \\
&G\ar[rd]&&\widehat G\ar[ld]&\\
&&\ast&&
}}
\end{equation}
of trivial bundles, where the top isomorphism $\pi$
is given by the pairing $\langle . , . \rangle:G\times\widehat G\to\U(1)$, i.e.
\begin{equation}\label{eqOfPo}
\pi(g,\chi,z):=(g,\chi, \langle g,\chi\rangle z).
\end{equation}
\end{defi}

A general Pontrjagin duality triple is a diagram which has 
the trivial Pontrjagin duality triple over the point as its local model.

\begin{defi}\label{defitDualityTriple}
A {\bf Pontrjagin duality triple}  over  $B$ 
is a commutative diagram 
\begin{equation}\label{diagDefiTdual}
\xymatrix{
&F\times_B \widehat E\ar[rd]\ar[ld]\ar[rr]^\kappa&&E\times_B \widehat F\ar[rd]\ar[ld]&\\
F\ar[rd]&&E\times_B \widehat E\ar[rd]\ar[ld]&&\widehat F\ar[ld]\\
&E\ar[rd]&&\widehat E\ar[ld]&\\
&&B&&
},
\end{equation}
where  $F\to E\to B$ is a   pair,
$\widehat F\to \widehat E\to B$ is a dual pair
and $\kappa$ is an isomorphism of $\U(1)$-bundles 
such that there exists atlases of the pair an the
dual pair over the same open cover $\{U_i\}$ of
$B$ such that for each point $b\in B$ there are 
charts containing $b$ such that the restrictions 
$$
\xymatrix{
F|_{E|_b}\ar[r]\ar[d]& G\times \U(1)\ar[d]&
\widehat F|_{\widehat E|_b}\ar[r]\ar[d]& \widehat G\times \U(1)\ar[d]\\
E|_b\ar[r]\ar[d]& G\ar[d]& 
\widehat E|_b\ar[r]\ar[d]&  \widehat G\ar[d]\\
b\ar[r]&\ast&b\ar[r]&\ast
}
$$
of these charts induce the trivial Pontrjagin 
duality diagram (\ref{diagDefiTrivTdual}) 
from the restriction (\ref{diagDefiTdual})$|_b$.

A {\bf chart} of a Pontrjagin duality triple
is the datum of two charts with the  property
spelled out above.
An {\bf atlas} of a Pontrjagin duality triple
is a collection of charts  covering $B$.  
\end{defi}

Let  $\hat p:\widehat E\to B$ be a $\widehat G$-bundle,
and let $\widehat E^{\rm op}$ be the opposite $\widehat G$-bundle
which is $\widehat E$ as a space but with $\widehat G$-action: 
$(\hat e,\chi)\mapsto \hat e\cdot (-\chi)$.
The simplest possible non-trivial example of a Pontrjagin 
duality triple is 
\begin{eqnarray}\label{DiagTheFirstExaOfAPDT}
\xymatrix{
&F_{\widehat E^{\rm op}}\times_B \widehat E\ar[rd]\ar[ld]\ar[rr]^{\kappa_{\widehat E}}&
&\ar[rd]\ar[ld](B\times G)&\hspace{-1,6cm} \times_B \widehat E\times \U(1)\\
F_{\widehat E^{\rm op}}\ar[rd]&&(B\times G)\times_B \widehat E\ar[rd]\ar[ld]&&\widehat E\times \U(1) \ar[ld]\\
&B\times G\ar[rd]&&\widehat E\ar[ld]_{\hat p}&\\
&&B&&
},
\end{eqnarray}
where 
$$
	\kappa_{\widehat E}( [\hat e,g,t],\hat e\cdot\chi):= 
	(\hat p(\hat e),g, \hat e\cdot\chi ,\langle g,\chi\rangle t)
$$
and the obvious maps elsewhere.
Note that it is necessary to deal with $F_{\widehat E^{\rm op}}$ instead of 
$F_{\widehat E}$ to make the top isomorphism $\kappa_{\widehat E}$
well-defined with the correct local structure (\ref{eqOfPo}).

Let $\{U_i\}$ be an atlas of a general Pontrjagin duality triple (\ref{diagDefiTdual}).
We have transition functions 
$g_{ji}:U_{ij}\to G,\zeta_{ji}:U_{ij}\times G\to\U(1)$ for the pair and
$\hat g_{ji}:U_{ij}\to\widehat G,\hat \zeta_{ji}:U_{ij}\times\widehat G\to\U(1)$
for the dual pair.
As  over each point the isomorphism $\kappa$ reduces
to $\pi$ from (\ref{diagDefiTrivTdual}),
the relation
\begin{eqnarray}
\hat \zeta_{ji}(u,\chi) \langle g,\chi\rangle = 
\langle g+g_{ji}(u),\chi+\hat\chi_{ji}(u)\rangle \zeta_{ji}(u,g)
\end{eqnarray}
holds.
By putting $\chi=0$ or $g=0$ we find that
\begin{eqnarray}
\zeta_{ji}(u,g) &=& \langle g,\hat \chi_{ji}(u)\rangle^{-1} \ s_{ji}(u),\\
\hat\zeta_{ji}(u,\chi)&=&\langle g_{ji}(u),\hat \chi\rangle\ 
\langle g_{ji}(u),\hat \chi_{ji}(u)\rangle\ 
s_{ji}(u),\nonumber
\end{eqnarray}
where $s_{ji}(u):=\zeta_{ji}(u,0)$.
Thus, if we compare with (\ref{EqDasGibtNenTorsorAlter}),
we see that Pontrjagin duality triples 
contain the same $\check C$ech theoretic
amount of data as module pairs.
In the following we make this correspondence more
precise by constructing an explicit equivalence 
of gerbes between the gerbe ${\rm Mod}_{E,\widehat E}$ of module pairs
of Proposition \ref{PropDieModGerbe} and the gerbe of Pontrjagin duality
triples we introduce next.

Let $B$ be a space and let
$E\to B$ and $\widehat E\to B$ be a $G$-bundle
and a $\widehat G$-bundle, respectively.
Let us denote by ${\rm Pon}(E,\widehat E)$ the groupoid
\label{PageOfPonTrip}
of Pontrjagin duality triples with $E\to B$ and $\widehat E\to B$ fixed,
i.e. a morphism 
$$
\resizebox*{5cm}{2,5cm}{$
\left(
\begin{array}{c}
\xymatrix{
&F\times_B \widehat E\ar[rd]\ar[ld]\ar[rr]^\kappa&&E\times_B \widehat F\ar[rd]\ar[ld]&\\
F\ar[rd]&&E\times_B \widehat E\ar[rd]\ar[ld]&&\widehat F\ar[ld]\\
&E\ar[rd]&&\widehat E\ar[ld]&\\
&&B&&
}
\end{array}
\right)$}
\to
\resizebox*{5cm}{2,5cm}{$
\left(
\begin{array}{c}
\xymatrix{
&F'\times_B \widehat E\ar[rd]\ar[ld]\ar[rr]^{\kappa'}&&E\times_B \widehat F'\ar[rd]\ar[ld]&\\
F'\ar[rd]&&E\times_B \widehat E\ar[rd]\ar[ld]&&\widehat F'\ar[ld]\\
&E\ar[rd]&&\widehat E\ar[ld]&\\
&&B&&
}\end{array}
\right)$
}
$$
between two Pontrjagin duality triples with fixed $E$ and $\widehat E$
consists of two bundle morphisms
\begin{eqnarray*}
\xymatrix{
F\ar[r]^\cong\ar[d]& F'\ar[d] &&\widehat F\ar[r]^\cong\ar[d]& \widehat F'\ar[d]\\
E\ar[r]^=& E&& \widehat E\ar[r]^=& \widehat E
}
\end{eqnarray*}
such that the induced diagram 
\begin{eqnarray*}
\xymatrix{
F\times_B {\widehat E}\ar[d]\ar[r]^\kappa & E\times_B F\ar[d]\\
F'\times_B \widehat E \ar[r]^{\kappa'}& E\times_B F'
}
\end{eqnarray*}
commutes.
For each open $U\subset B$ let 
${\rm Pon}_{E,\widehat E}(U):={\rm Pon}(E|_U,\widehat E|_U)$ be
the groupoid of Pontrjagin duality triples over $E|_U\to U$ and $\widehat E|_U\to U$.

\begin{prop}
The assignment
$$
{\rm Pon}_{E,\widehat E}:U\mapsto {\rm Pon}(E|_U,\widehat E|_U)
$$ 
is a $\U(1)$-banded gerbe.
\end{prop}

\begin{pf}
The proof is straight forward and analogue to the proof of Proposition \ref{PropDieModGerbe}.
\end{pf}

As in the proof of Proposition \ref{PropDieModGerbe}, we can also determine 
the class of the gerbe ${\rm Pon}_{E,\widehat E}$, but the result of next section  
will do this job equally well. Namely, we construct an explicit equivalence of gerbes
${\rm Mod}_{E,\widehat E}\cong{\rm Pon}_{E,\widehat E}.$

\section{Module pairs vs. Pontrjagin duality triples}

\noindent
Let $\widehat E\to B$ be a $\widehat G$-bundle, and  
let $F_{\widehat E}\to B\times G\to B$ be the ring pair as defined
in (\ref{EqGetMyFirstRingStr}), i.e. $F_{\widehat E}:=
\widehat E\times_{\widehat G}( G\times \U(1))$.
There is a canonical map $\iota:\widehat E\times G\to F_{\widehat E}$
given by  $\iota(\hat e,g):=[\hat e,g,1]$.
Let $F\stackrel{q}{\to} E\to B$ be a module pair over $F_{\widehat E}\to B\times G\to B$
with action map $\varrho:F\times_B F_{\widehat E}\to F$.
By use of the maps $\rho$ and $\iota$ we define a $G$-action on
the space $F\times_B\widehat E$ by
\begin{eqnarray}\label{EqGactionOnEtimesF}
(F\times_B\widehat E)\times G &\to& F\times_B\widehat E.\\
(x,\hat e,g)&\mapsto& (\varrho(x,\iota(\hat e,g) ), \hat e)\nonumber
\end{eqnarray}
The corresponding quotient 
$
\widehat F_\varrho:= \big(F\times_B\widehat E\big)\big/{G}
$
comes along with two natural maps:
\begin{eqnarray*}
\begin{array}{rccccccrcc}
\hat q :\!\!\!\!& \widehat F_\varrho&\to&\widehat E&&{\rm and}&&
\kappa_\varrho:\ F\times_B \widehat E&\to&E\times_B\widehat F_\varrho.\\
&{[x,\hat e]}&\mapsto &\hat e&&&&
(x,\hat e)&\mapsto& (q(x), [x,\hat e])
\end{array}
\end{eqnarray*}
The map $\hat q$ is a principal $\U(1)$-bundle
with action induced by the principal action of 
$\U(1)$ on $F$.

\begin{prop}\label{PropReturnToSender}
The diagram 
\begin{eqnarray}\label{DiagPDTinPropos}
\xymatrix{
&F\times_B \widehat E^{\rm op}\ar[rd]\ar[ld]\ar[rr]^{\kappa_\varrho}&
&E\times_B \widehat F_\varrho\ar[rd]\ar[ld]&\\
F\ar[rd]^q&&E\times_B\widehat E^{\rm op}\ar[rd]\ar[ld]&
&\widehat F_\varrho\ar[ld]_{\hat q}\\
&E\ar[rd]&&\widehat E^{\rm op}\ar[ld]&\\
&&B&&
}
\end{eqnarray}
is a Pontrjagin duality triple.
Here $\widehat E^{\rm op}$ the space
$\widehat E$ equipped with the opposite $\widehat G$-action
$(\hat e,\chi)\mapsto \hat e\cdot(-\chi)$.
\end{prop}

\begin{pf}
Let $U_i$ be a module atlas, so there are  trivialisations 
$\widehat E|_{U_i}\cong U_i\times \widehat G$ and
$F|_{E|_{U_i}}\cong U_i\times G\times \U(1)$ such that
over each point of $U_i$ the action map induced from the 
action map $\varrho$ takes the form
$+\times\cdot:\U(1)\times G \times \U(1)\times G\to \U(1)\times G$.
These induce trivialisations 
$$
\widehat F_\varrho|_{\widehat E|_{U_i}}\cong 
\big((U_i\times G\times \U(1))\times_{U_i} (U_i\times \widehat G)\big)\big/ G,
$$
where the quotient of the right hand side is by the induced action:
$$
((u,g,t),(u,\chi) )\cdot h :=((u,g+h,\langle h,\chi\rangle t),(u,\chi)).
$$
We can identify further
\begin{eqnarray*}
\big((U_i\times G\times \U(1))\times_{U_i} (U_i\times \widehat G)\big)\big/ G
&\cong&
U_i\times \widehat G\times\U(1)\\
\ [(u,g,t),(u,\chi)]&\mapsto &(u,\chi,\langle g,\chi\rangle^{-1} t)
\end{eqnarray*}

From this one can deduce that, firstly,   
$\widehat F_\varrho\to\widehat E\to B$ is a dual pair, and that, secondly,
the isomorphism $\kappa_\varrho$ induces 
\begin{eqnarray*}
\xymatrix{
(F\times_B\widehat E)|_{U_i}\ar[rr]^{\kappa_\varrho|_{U_i}}\ar[d]&
&(E\times_B \widehat F_\varrho)|_{U_i}\ar[d]\\
\big(U_i\times G\times \U(1)\big)\times_{U_i}\big(U_i\times\widehat G\big)
\ar[rr]&&
\big(U_i\times G\big)\times_{U_i} \big(U_i\times \widehat G\times\U(1)\big)
}
\end{eqnarray*}
\vspace{-0,6cm}
$$
\quad\ 
\xymatrix{
\big( (u,g,t),(u,\chi)\big)\ar@{|->}[rr]&& \big( (u,g),(u,\chi,\langle g,\chi\rangle^{-1} t)\big)
}
$$
where the vertical arrows are the chosen trivialisations.
It follows  that $\kappa_\varrho:F\times_B \widehat E\to E\times_B\widehat F_\varrho$
does not have the correct local structure of (\ref{eqOfPo}), due to the
power of $-1$ in the above diagram.
However, the minus sign can be absorbed by a modification of the trivialisations 
by  the map $\chi\mapsto -\chi$ (which is not a $\widehat G$-bundle morphism).
This means nothing but dealing with  $\widehat E^{\rm op}$ 
instead of $\widehat E$.
This proves the proposition.
\end{pf}

The above construction is completly natural, therfore  
it defines a morphism of gerbes
\begin{eqnarray}\label{EqFromModToPon}
	{\rm Mod}_{E,\widehat E}\to{\rm Pon}_{E,\widehat E^{\rm op}}
	\textrm{  or likewise  }
	{\rm Mod}_{E,\widehat E^{\rm op}}\to{\rm Pon}_{E,\widehat E}.
\end{eqnarray}
An  example of the above construction is the following.
Let $\widehat E\to B$ be a $\widehat G$-bundle.
Consider the opposite $\widehat E^{\rm op}$ and note that 
$(\widehat E^{\rm op})^{\rm op}=\widehat E$.
Take the ring pair $F_{\widehat E^{\rm op}}\to B\times G\to B$ 
as a module pair over itself. Then the above construction 
yields a Pontrjagin duality triple 
\begin{eqnarray*}
	\xymatrix{
	&F_{\widehat E^{\rm op}}\times_B \widehat E\ar[rd]\ar[ld]\ar[rr]^{\kappa}&
	&(B\times G)\ar[rd]\ar[ld]&\hspace{-1,0cm}
	\times_B (F_{\widehat E^{\rm op}}\times_B\widehat E)/G&\\
	F_{\widehat E^{\rm op}}\ar[rd]&&(B\times G)\times_B \widehat E\ar[rd]\ar[ld]&
	&(F_{\widehat E^{\rm op}} \ar[ld]&\hspace{-1,6cm}\times_B\widehat E)/G\\
	&B\times G\ar[rd]&&\widehat E\ar[ld]&&\\
	&&B&&&
	}
\end{eqnarray*}
which is naturally isomorphic to (\ref{DiagTheFirstExaOfAPDT})
by the $\U(1)$-bundle isomorphism 
\begin{eqnarray*}
	(F_{\widehat E^{\rm op}}\times_B\widehat E)/G&\cong& \widehat E\times\U(1).\\
	\ [[\hat e,g,t],\hat e\cdot \chi]&\mapsto& (\hat e\cdot\chi, \langle g,\chi\rangle t)
\end{eqnarray*}

Let us start with the construction of a morphism in opposite direction 
of (\ref{EqFromModToPon}). Consider a Pontrjagin duality triple
$$
\xymatrix{
&F\times_B \widehat E\ar[rd]\ar[ld]_{{\rm pr}_F}\ar[rr]^\kappa&&E\times_B \widehat F\ar[rd]\ar[ld]&\\
F\ar[rd]&&E\times_B \widehat E\ar[rd]\ar[ld]&&\widehat F\ar[ld]\\
&E\ar[rd]&&\widehat E\ar[ld]&\\
&&B&&
}.
$$
The $\widehat G$-bunde $\widehat E$ defines the ring pair 
$F_{\widehat E^{\rm op}}\to B\times G\to B$,
i.e. $F_{\widehat E^{\rm op}}:= \widehat E^{\rm op}\times_{\widehat G} (G\times\U(1))$,
and we define 
\begin{eqnarray}\label{EqWerhettgeadcachtdat}
\varrho_\kappa:F\times_B F_{\widehat E^{\rm op}}&\to &F.\\
(x,[\hat e,g,t])&\mapsto& {\rm pr}_F( \kappa^{-1}(\kappa(x,\hat e)\cdot g)) \cdot t\nonumber
\end{eqnarray}
This is well-defined, because  $\kappa$ is
locally of the form (\ref{eqOfPo}) and  therefore
satisfies
\begin{eqnarray*}
\kappa^{-1}(\kappa(x,\hat e\cdot \chi)\cdot g)=
\big(\kappa^{-1}(\kappa(x,\hat e)\cdot g)\cdot\chi\big)\cdot \langle g,\chi\rangle^{-1}.
\end{eqnarray*}
Here the action of $G$ on $E\times_B \widehat F$ and the actions 
of  $\widehat G$ and $\U(1)$ on $F\times_B \widehat E$ are the 
obvious ones.

\begin{prop}
Via the map $\varrho_\kappa$ the pair $F\to E\to B$ 
is a module pair over the ring pair $F_{\widehat E^{\rm op}}\to B\times G\to B$.
\end{prop}

\begin{pf}
It is immediate that $\varrho_\kappa$ satisfies 
the necessary local condition.
\end{pf}

As the construction of $\varrho_\kappa$ is natural, 
we obtain a morphism of gerbes
\begin{eqnarray}\label{EqFromPonToMod}
{\rm Pon}_{E,\widehat E}\to{\rm Mod}_{E,\widehat E^{\rm op}}.
\end{eqnarray}

\begin{prop}\label{PropEqivalenceOfGerbes}
$
(\ref{EqFromModToPon}):{\rm Mod}_{E,\widehat E^{\rm op}}
\rightleftarrows
{\rm Pon}_{E,\widehat E}:(\ref{EqFromPonToMod})
$
is an equivalence of gerbes.
\end{prop}

\begin{pf}
Consider the composition $(\ref{EqFromPonToMod})\circ(\ref{EqFromModToPon})$
applied to a module pair $F\to E\to B$ with action map $\varrho:F\times_B F_{\widehat E}\to F$.
This gives a new module pair, where the underlying pairs are unchanged, but with 
a new action map $\varrho_{\kappa_\varrho}: F\times F_{\widehat E^{\rm op}}\to F$
as obtained from the top isomorphism 
$\kappa_\varrho$. Note that $\kappa_\varrho$ is $G$-equivariant with respect to 
the $G$-action  (\ref{EqGactionOnEtimesF}) on $F\times_B \widehat E$
and the obvious $G$-action on $E\times_B \widehat F.$
Therefore we just can compute $\varrho_{\kappa_\varrho}$:
\begin{eqnarray*}
\varrho_{\kappa_\varrho}(x,[\hat e,g,t]) 
&=&{\rm pr}_F (\kappa_\varrho^{-1}(\kappa_\varrho(x,\hat e)\cdot g))\cdot t\\
&=&{\rm pr}_F (\varrho(x,\iota(\hat e,g)),\hat e) \cdot t\\
&=&\varrho(x,\iota(\hat e,g)) \cdot t\\
&=&\varrho(x,[\hat e,g,t]).
\end{eqnarray*}
Thus, the composition $(\ref{EqFromPonToMod})\circ(\ref{EqFromModToPon})$
is the identity on ${\rm Mod}_{E,\widehat E^{\rm op}}$.

Consider the composition $(\ref{EqFromModToPon})\circ(\ref{EqFromPonToMod})$
applied to a Pontrjagin duality triple with top isomorphism
$\kappa:F\times_B\widehat E\to E\times \widehat F$.
Define $\varrho_\kappa$,
$\widehat F_{\varrho_\kappa}$ and
$\kappa_{\varrho_\kappa}:F\times_B\widehat E\to E\times \widehat F_{\varrho_\kappa}$
as above. Note that $\kappa:F\times_B\widehat E\to E\times_B \widehat F$
is $G$-equivariant with respect to the $G$-action (\ref{EqGactionOnEtimesF})
induced by $\varrho_\kappa$ and the obvious $G$-action on $E\times_B \widehat F$.
Thus, there are well-defined bundle isomorphims 
\begin{eqnarray*}
\xymatrix{
F\ar[r]^=\ar[d]& F\ar[d] &&
\widehat F_{\varrho_\kappa}\ar[rrr]^{[x,\hat e]\mapsto{\rm pr}_{\widehat F}(\kappa(x,\hat e))}\ar[d]
&&& \widehat F\ar[d]   \\
E\ar[r]^=& E&& \widehat E\ar[rrr]^=&&& \widehat E
},
\end{eqnarray*}
and the induced diagram  
\begin{eqnarray*}
\xymatrix{
F\times_B {\widehat E}\ar[d]^=\ar[r]^{\kappa_{\varrho_\kappa}} & E\times_B F_{\varrho_\kappa}\ar[d]^\cong\\
F\times_B \widehat E \ar[r]^{\kappa}& E\times_B F
}
\end{eqnarray*}
commutes. This construction is a natural isomorphism of 
Pontrjagin duality triples, therefore the composition
$(\ref{EqFromModToPon})\circ(\ref{EqFromPonToMod})$
is naturally isomorphic to the identity transformation on ${\rm Pon}_{E,\widehat E}$.
\end{pf}

\begin{cor}\label{CorTheClassOfPon}
	Let  $E\to B$ and  $\widehat E\to B$ be a $G$-bundle 
	and a  $\widehat G$-bundle, respectively.
	\begin{itemize}
	\item[$(i)$]
		The class of the gerbe ${\rm Pon}_{E,\widehat E}$ is given 
		by $-[E]\cup[\widehat E]\in 
		\check H^2(B,\underline{\U(1)})$.
	\item[$(ii)$]
		If $[E]\cup[\widehat E]=0$, then the set of isomorphism classes of 
		the groupoid ${\rm Pon}(E,\widehat E)$ is a
		$\check H^1(B,\underline{\U(1)})$-torsor.
	\end{itemize}
\end{cor}

\begin{pf}
	Proposition \ref{PropEqivalenceOfGerbes}
	and Proposition \ref{PropDieModGerbe} together with the fact
	$[\widehat E^{\rm op}]=-[\widehat E]$ imply $(i)$.
	Proposition \ref{PropEqivalenceOfGerbes}
	and Proposition \ref {PropIsoClassesOfMod}
	imply $(ii)$.
\end{pf}

\begin{exa}
Let $B:=\Sigma$ be a two-dimensional, connected, closed manifold,
then $\check H^2(\Sigma,\underline{\U(1)})\cong H^3(\Sigma,\Z)=0$,
so all choices of bundles $E\to \Sigma\leftarrow\widehat E$
can be extended to a Pontrjagin dualitiy triple.

If $\Sigma$ is non-orientalbe, then 
also $\check H^1(\Sigma,\underline{\U(1)})\cong H^2(\Sigma,\Z)=0$,
so up to isomorphism in ${\rm Pon}(E,\widehat E)$ this 
Pontrjagin duality triple is unique.
Otherwise, if $\Sigma$ is orientable, then
$\check H^1(\Sigma,\underline{\U(1)})\cong H^2(\Sigma,\Z)\cong\Z$,
so there are infinitely many different isomorphism classes
in ${\rm Pon}(E,\widehat E)$.
\end{exa}

\section{Extensions to Pontrjagin duality triples}
\label{SecExtensions}

\begin{defi}
Let $X$ be a commutative diagram of topological spaces. 
An {\bf extension} of $X$ to a Pontrjagin duality triple
is a Pontrjagin duality triple which contains $X$ as a
sub-diagram.
\end{defi}

Let $E\to B$ be a $G$-bundle, and $\widehat E\to B$
be a $\widehat G$-bundle.
The first part of Corollary \ref{CorTheClassOfPon} answers 
the question whether or not the diagram 
\begin{eqnarray}\label{DiagTheFirstExtProbl}
\resizebox*{!}{3cm}{$
\xymatrix{
&E\times_B\widehat E\ar[dr]\ar[dl]&\\
E\ar[dr]&&\widehat E\ar[dl]\\
&B&
}$}
\end{eqnarray}
admits an extension. Namely, 
a global object of the gerbe ${\rm Pon}_{E,\widehat E}$ exists,
 i.e. it exists an object in ${\rm Pon}(E,\widehat E)=
 {\rm Pon}_{E,\widehat E}(B),$
if and only if the class of the gerbe is trivial. 
Thus, diagram (\ref{DiagTheFirstExtProbl}) can be extended
to a Pontrjagin duality triple if and only if $[E]\cup[\widehat E]=0$.
The second part of Corollary \ref{CorTheClassOfPon} 
states that such extensions are not unique if they exist, but
(up to isomorphism) we know exactly about this ambiguity
which only depends on the topology of $B$ 
(and not on $G$ or $\widehat G$).

We wish to understand another extension problem. Namely,
if given a pair $ F\to E\to B$, we  want to understand 
the existence and uniqueness problem of extensions of the
diagram 
\begin{eqnarray}\label{DiaTheExtProbOfPairs}
\resizebox*{!}{3cm}{$
\begin{array}{c}
\xymatrix{
F\ar[rd]&&\\
&E\ar[rd]&\\
&&B
}
\end{array}
$}.
\end{eqnarray}
To manage this we take a closer look at the category ${\rm Pon}(B)$
of (all) Pontrjagin duality triples over a space $B$.
A morphisms in this category 
\begin{eqnarray}\label{DiagCatOfAllPDT}
(a,\hat a):\resizebox*{5cm}{2,7cm}{$
\left(
\begin{array}{c}
\xymatrix{
&F\times_B \widehat E\ar[rd]\ar[ld]\ar[rr]^\kappa&&E\times_B \widehat F\ar[rd]\ar[ld]&\\
F\ar[rd]&&E\times_B \widehat E\ar[rd]\ar[ld]&&\widehat F\ar[ld]\\
&E\ar[rd]&&\widehat E\ar[ld]&\\
&&B&&
}
\end{array}
\right)$}
\to
\resizebox*{5cm}{2,7cm}{$
\left(
\begin{array}{c}
\xymatrix{
&F\boldsymbol{'}\times_B \widehat E\boldsymbol{'}
\ar[rd]\ar[ld]\ar[rr]^{\kappa\boldsymbol{'}}&
&E\boldsymbol{'}\times_B \widehat F\boldsymbol{'}\ar[rd]\ar[ld]&\\
F\boldsymbol{'}\ar[rd]&&E\boldsymbol{'}\times_B \widehat E\boldsymbol{'}\ar[rd]\ar[ld]&&\widehat F\boldsymbol{'}\ar[ld]\\
&E\boldsymbol{'}\ar[rd]&&\widehat E\boldsymbol{'}\ar[ld]&\\
&&B&&
}\end{array}
\right)$
}
\end{eqnarray}
consists of a morphism 
of pairs  $a:(F\to E\to B) \to(F'\to E'\to B)$ over $B$ and of 
a morphism of dual pairs  $\hat a:(\widehat F\to \widehat E\to B)
\to(\widehat F'\to\widehat E'\to B)$ over $B$ such that the 
induced diagram
\begin{eqnarray*}
\xymatrix{
F\times_B {\widehat E}\ar[d]\ar[r]^\kappa & E\times_B F\ar[d]\\
F'\times_B \widehat E' \ar[r]^{\kappa'}& E'\times_B F'
}
\end{eqnarray*}
commutes.
Recall from section \ref{SecPairs} that the
automorphism group of the trivial pair over 
the point is the semi-direct product 
${\rm A_{\rm Par}}=G\ltimes C(G,\U(1))$.
The automorphism group of the trivial dual pair over the
point is $ \widehat{\rm  A}_{\rm Par}:=\widehat G\ltimes C(\widehat G,\U(1))$.
These two groups contain two isomorphic subgroups
\begin{eqnarray*}
{\rm A_{\rm Par}}\supset G\ltimes (\U(1)\times \widehat G)\stackrel{\phi}{\cong}
 \widehat G\ltimes (\U(1)\times G)\subset \widehat{\rm A}_{\rm Par}.
\end{eqnarray*}
Here we take the isomorphism  $\phi:(g,t,\chi)\mapsto (-\chi, t\langle g,-\chi\rangle,g)$,
and $\U(1)\times\widehat G$ is the subgroup of $C(G,\U(1))$ consisting of those $f$ such that
$f(h)=t\langle h,\chi\rangle$ for some $t\in\U(1)$ and $\chi\in\widehat G$. 
Similarly, the inclusion $\U(1)\times G\subset C(\widehat G,\U(1))$
is understood by the identification $G\cong \widehat{\widehat G}$.

\begin{prop}
	The automorphism group of the 
	trivial Pontrjagin duality 
	triple over the point is 
	\begin{eqnarray*}
		{\rm A_{\rm Pon}}:= \{ (a,\phi(a))| a\in G\ltimes (\U(1)\times \widehat G)\}.
	\end{eqnarray*}
\end{prop}

\begin{pf}
	Given $(g,f)\in {\rm A_{\rm Par}}$ and $(\chi,F)\in \widehat{\rm A}_{\rm Par}$, 
	then it is easy to see that 
	$$
	\xymatrix{
	G\times \U(1)\times {\widehat G}\ar[d]^{(g,f)\times\chi}\ar[r]^\pi & 
	G\times \widehat G\times	\U(1)\ar[d]^{g\times(\chi,F)}\\
	G\times\U(1)\times \widehat G \ar[r]^{\pi}& G\times \widehat G\times \U(1)
	}
	$$
	commutes if and only if $f(h)= f(0)\langle h,\chi\rangle^{-1}$
	and $F(\psi)=f(0)\langle g,\chi\rangle\langle g,\psi\rangle$.
	This encodes precisely the isomorphism $\phi$.
\end{pf} 

To an ${\rm A_{\rm Pon}}$-principal bundle $P\to B$ we can associate 
a Pontrjagin duality triple 
$$
P\times_{\rm A_{\rm Pon}} 
\boldsymbol{\resizebox*{4cm}{2,3cm}{$
\left(
\begin{array}{c}
\xymatrix{
&G\times \widehat G\times\U(1)\ar[rd]\ar[ld]\ar[rr]^{\pi}&
&G\times \widehat G\times\U(1)\ar[rd]\ar[ld]&\\
G\times \U(1)\ar[rd]&&G\times \widehat G\ar[rd]\ar[ld]&&
\widehat G\times\U(1)\ar[ld]\\
&G\ar[rd]&&\widehat G\ar[ld]&\\
&&\ast&&
}\end{array}
\right)$
}}
$$
over $B$ which defines a functor from the category of ${\rm A_{\rm Pon}}$-principal
fibre bundles over $B$ to the category ${\rm Pon}(B)$ of Pontrjagin duality triples over $B$.

\begin{prop}\label{PropEqOfCatsPonAndAutoOfPon}
	The functor
	\begin{eqnarray*}
		(P\to B)\mapsto 
		P\times_{\rm A_{\rm Pon}} 
		\boldsymbol{
		\resizebox*{4cm}{2,3cm}{$
		\left(
		\begin{array}{c}	
		\xymatrix{
		&G\times \widehat G\times\U(1)\ar[rd]\ar[ld]\ar[rr]^{\pi}&
		&G\times \widehat G\times\U(1)\ar[rd]\ar[ld]&\\
		G\times \U(1)\ar[rd]&&G\times \widehat G\ar[rd]\ar[ld]&&
		\widehat G\times\U(1)\ar[ld]\\
		&G\ar[rd]&&\widehat G\ar[ld]&\\
		&&\ast&&
		}\end{array}
		\right)$
		}}
	\end{eqnarray*}
	is an equivalence of categories.
\end{prop}

\begin{pf}
	The proof is analogous to the proofs of 
	Proposition \ref{PropEqOfCatsPairsAndBunAFull}
	and Proposition \ref{PropMyFirstEqOfCats}.
	We sketch how to define an inverse functor
	up to equivalence.
	Let 
	\begin{eqnarray*}
	\resizebox*{!}{2,7cm}{$
	\begin{array}{c}
	\xymatrix{
	&F\times_B \widehat E\ar[rd]\ar[ld]\ar[rr]^\kappa&
	&E\times_B \widehat F\ar[rd]\ar[ld]&\\
	F\ar[rd]_q&&E\times_B \widehat E\ar[rd]\ar[ld]&&\widehat F\ar[ld]^{\hat q}\\
	&E\ar[rd]&&\widehat E\ar[ld]&\\
	&&B&&
	}
	\end{array}
	$}
	\end{eqnarray*}
	be a Pontrjagin duality triple. For a $b\in B$ we 
	define $P_b$ to be the set of all tuples $(e,s,\hat e,\hat s)$
	with $e\in E|_b, \hat e\in \widehat E|_b$ and
	$s\subset F|_{E|_b}$, $\hat s\subset \widehat F|_{\widehat E|_b}$
	such that 
	$q|_s: s\to E|_b$ and $\hat q|_{\hat s}:\hat s\to \widehat E|_b$
	are homoeomorphisms, and such that $\kappa$ restricts to 
	homoeomorphisms
	$$
	\{s(e)\}\times \widehat E|_b\cong \{e\}\times\hat s,\qquad 
	s\times \{\hat e\}\cong E\times \{\hat s(\hat e)\},
	$$
	where $s(e):=(q|_{s})^{-1}(e)$ and $\hat s(\hat e):=(\hat q|_{\hat s})^{-1}(\hat e)$.
	$P_b$ is a ${\rm A_{Pon}}$-torsor subject to the action
	\begin{eqnarray*}
		P_{b}\times {\rm BA_{Pon}}&\to& P_{b},\\
		\big((e,s,\hat e,\hat s),(a,\phi(a))\big)&\mapsto &\big((e,s)
		\diamond\hspace{-0.24cm}\cdot\ a,		(\hat e,\hat s)
		\hat \diamond\hspace{-0.24cm}\cdot\ \phi(a)\big),
	\end{eqnarray*}
	where $\diamond\hspace{-0.165cm}\cdot$ is as in
	(\ref{EqNetteActionImZugNachWien}), and 
	$\hat\diamond\hspace{-0.165cm}\cdot$ is defined by the same 
	formula but with $G$ and $\widehat G$ exchanged. 
	In fact, this action is 
	well-defined, free and transitive. The local trivialisations 
	of the triple induce on $\coprod_{b\in B} P_b$ a topology such that 
	the projection $\coprod_{b\in B} P_b\to B$ becomes a principal
	${\rm A_{Pon}}$-fibre bundle.
	This way we obtain a functor from Pontrjagin duality triples 
	to ${\rm A_{Pon}}$-bundles which is up to equivalence inverse 
	to the functor in the proposition.
\end{pf}

Let us denote by ${\rm Pon^{full}}(F,E)$ the full subcategory of ${\rm Pon}(B)$
whose objects are the extensions of the pair $F\to E\to B$. By 
${\rm Pon}(F,E)$ we denote the proper subcategory of ${\rm Pon}(B)$
whose objects are the extensions of $F\to E\to B$ and whose 
morphisms (\ref{DiagCatOfAllPDT})
are only those of the form $(a,\hat a)=( {\rm id},\hat a),$
i.e. only those which are the identity on the underlying pair.
The task of understanding the extension problem of 
a pair $F\to E\to B$ is to understand the 
the categories ${\rm Pon}^{\rm full}(F,E)$ and ${\rm Pon}(F,E)$ or, at least, to understand the structure 
of their isomorphism classes $[{\rm Pon}^{\rm full}(F,E)]$ and $[{\rm Pon}(F,E)]$.
We will see next that this is equivalent to understand corresponding categories of 
${\rm A_{Pon}}$-reductions of the ${\rm A_{Par}}$-principal bundle
given by the pair (Proposition \ref{PropEqOfCatsPairsAndBunAFull}).

Let $P\to B$ be a principal ${\rm A_{Par}}$-bundle. The category of reductions
${\rm Red}(P)$ of $P$ to ${\rm A_{Pon}}$-bundles has as objects
commutative diagrams
\begin{eqnarray}\label{EqThisIsARedNotGreen}
\xymatrix{
P_{r}\ar[r]\ar[d]&P\ar[d]\\
B\ar[r]^=&B
},
\end{eqnarray}
where $P_{r}\to B$ is a ${\rm A_{Pon}}$-principal bundle and
$P_{r}\to P$ is an ${\rm A_{Pon}}$-equivariant map with closed
image.
A morphism in this category is a commutative diagram 
$$
\xymatrix{
P_{r}\ar[rd]\ar[dd]\ar[rrrd]&&&\\
&P_{s}\ar[rr]\ar[dd]&&P\ar[dd]\\
B\ar[rd]_{=}\ar[drrr]|!{[r];[dr]}\hole^=&&&\\
&B\ar[rr]^=&&B
}
$$
where $P_{r}\to P_{s}$ is an ${\rm A_{Pon}}$-bundle isomorphism over $B$.

\begin{prop}\label{PropEqOfCatsForRed}
	Let $P\to B$ be an ${\rm A_{Par}}$-bundle, and let 
	$F\to E\to B$ be the pair associated to $P$ by the functor of
	Proposition \ref{PropEqOfCatsPairsAndBunAFull}.
	Then there is an equivalence of categories 
	$$
	{\rm Red}(P)\simeq {\rm Pon}(F,E).
	$$
\end{prop}

\begin{pf}
	The proof is straight forward. From a reduction
	$P_{r}$ of $P$ one obtains an extension $X$ of 
	$F\to E\to B$ by associating the trivial
	Pontrjagin duality triple to $P_{r}$ and by observing
	that $F\to E\to B$ and the pair of $X$ are isomorphic
	by the ${\rm A_{Pon}}$-equivariant map $P_{r}\to P$.
	
	From an extension $X$ of $F\to E\to B$ one obtains a reduction
	of $P$ by the functor from ${\rm Pon}(B)$ to ${\rm A_{Pon}}$-bundles
	as constructed in the proof of Proposition \ref{PropEqOfCatsPonAndAutoOfPon}.
	
	Keeping track of the particular morphisms in the two categories 
	one obtains two functors  by this procedure which are inverses 
	of each other (up to natural equivalence).
\end{pf}

The isomorphism classes of the category of reductions 
are  a well-known object. If the quotient map ${\rm A_{Par}}\to {\rm A_{Par}}/{\rm A_{Pon}}$
has local sections (e.g. if G is compact), then by 
\cite[Theorem V.3.1]{Br} we have that the isomorphism classes $[{\rm Red}(P)]$
of reductions of $P$ are naturally bijective to the set
of sections $\Gamma(B,P/{\rm A_{Pon}})$ of 
$P/{\rm A_{Pon}} \cong P\times_{\rm A_{Par}}({\rm A_{Par}}/{\rm A_{Pon}})
\to B$ which is a fibre bundle with fibre ${\rm A_{Par}}/{\rm A_{Pon}}$, 
and we can identify the space of sections
with the ${\rm A_{Par}}$-equivariant maps $C(P, {\rm A_{Par}}/{\rm A_{Pon}})^{\rm A_{Par}}$.
However, if the quotient ${\rm A_{Par}}/{\rm A_{Pon}}$ behaves badly,
i.e. the quotient map does not have local sections, then 
the homotopty quotients
${\rm A_{Par}}\sslash{\rm A_{Pon}}:=({\rm A_{par}}\times {\rm EA_{Pon}})/{\rm A_{Pon}}$
and
$P\sslash{\rm A_{Pon}}:=(P\times {\rm EA_{Pon}})/{\rm A_{Pon}}
\cong P\times_{\rm A_{Par}} ({\rm A_{Par}}\sslash{\rm A_{Pon}})$
must be concerned. 
The technical issue we have to take care of here is 
that $P\to P/{\rm A_{Pon}}$ is only a principal bundle,
i.e. it has local sections, if ${\rm A_{Par}}\to {\rm A_{Par}}/{\rm A_{Pon}}$
has.

Generalising the proof of \cite[Theorem V.3.1]{Br}
we obtain the following classification result.

\begin{prop}\label{PropTheClassOfRedsAndExts}
	Let $P$ and $F\to E\to B$ be as above.
	\begin{itemize}
	\item[i)]
		There are natural bijections 
		$$
		[{\rm Pon}(F,E)]\cong [{\rm Red}(P)]\cong {\rm im}(\gamma_{*})\cong {\rm im}(\varepsilon_{*}),
		$$
		where 
		$\gamma_{*}:\Gamma(B,P\sslash{\rm A_{Pon}})\to\Gamma(B,P/{\rm A_{Pon}})$
		is the induced map of the canonical map $\gamma:P\sslash{\rm A_{Pon}}\to P/{\rm A_{Pon}}$,
		and 
		$\varepsilon_{*}:C(P, {\rm A_{Par}}\sslash{\rm A_{Pon}})^{\rm A_{Par}}
		\to C(P,{\rm A_{Par}}/{\rm A_{Pon}})^{\rm A_{Par}}$
		 is induced by $\varepsilon:{\rm A_{Par}}\sslash{\rm A_{Pon}}\to {\rm A_{Par}}/{\rm A_{Pon}}$.
	\item[ii)] If ${\rm A_{Par}}\to {\rm A_{Par}}/{\rm A_{Pon}}$ has local sections, then $\gamma_{*}$ and
		$\varepsilon_{*}$ are surjective, so
		$$
		[{\rm Pon}(F,E)]\cong [{\rm Red}(P)]\cong \Gamma(B,P/{\rm A_{Pon}})
		\cong C(P,{\rm A_{Par}}/{\rm A_{Pon}})^{\rm A_{Par}},
		$$
		
	\end{itemize}
\end{prop}

\begin{pf}
	$i)$ The first and the last bijection are immediate. 
	
	To construct a map
	from $[{\rm Red}(P)]$ to ${\rm im}(\gamma_{*})$ let $P_{r}\to P$ be a reduction
	of $P$, it induces a map $\sigma:B=P_{r}/{\rm A_{Pon}}\to P/{\rm A_{Pon}}$. 
	$\sigma$ is a section of $P/{\rm A_{Pon}}$ and does not vary inside the isomorphism
	class of the reduction.
	We show that $\sigma$ is in the image of $\gamma_{*}$.
	Choose any classifying bundle map $P_{r}\to {\rm EA_{Pon}}$,
	then we have an equivariant factorisation 
	$
	P_{r}\to P\times {\rm EA_{Pon}} \to P
	$
	of the reduction $P_{r}\to P$. It induces a factorisation  
	$
	B\to P\sslash {\rm A_{Pon}}\to P/{\rm A_{Pon}}
	$
	of $\sigma.$
	
	Given a section $\sigma\in{\rm im}(\gamma_{*})$ we define
	a reduction $P_{\sigma}$ to be the pullback in
	$$
	\xymatrix{
	P_{\sigma}\ar[r]\ar[d]& P\ar[d]\\
	B\ar[r]^-{\sigma} &P/{\rm A_{Pon}}
	}.
	$$
	Despite the fact that $P\to P/{\rm A_{Pon}}$ may not have local
	sections the pullback $P\times {\rm EA_{Pon}}$ 
	along $\gamma$ has (as indicated by the dashed arrow)
	$$
	\xymatrix{
	P\times{\rm EA_{Pon}}\ar[r]\ar[d]& P\ar[d]\\
	P\sslash{\rm A_{Pon}}\ar@{-->}@/_0.4cm/[u] \ar[r]^\gamma& P/{\rm A_{Pon}}
	},
	$$
	and as $\sigma$ has a factorisation over $\gamma$
	we find that $P_{\sigma}$ is also the pullback in 	
	$$
	\xymatrix{
	P_{\sigma}\ar[r]\ar[d]& P\times{\rm EA_{Pon}}\ar[r] \ar[d]&P\ar[d]\\
	B\ar@/_{0.8cm}/[rr]^\sigma\ar[r]\ar@{-->}@/_0.4cm/[u] 
	&P\sslash{\rm A_{Pon}}\ar@{-->}@/_0.4cm/[u]\ar[r]^\gamma&P/{\rm A_{Pon}}
	}.
	$$
	Therefore $P_{\sigma}\to B$ has local sections and is a principal ${\rm A_{Pon}}$-bundle.

	The two constructions  are easily seen to be inverses of each other. 
	
	$ii)$ In case $P\to P/{\rm A_{Pon}}$ is a pricipal ${\rm A_{Pon}}$-bundle
	it has local sections, and each section in $\Gamma(B,P/{\rm A_{Pon}})$ factors
	over $\gamma$.
\end{pf}

To classify the extensions of a pair $F\to E\to B$ up to isomorphism 
in the category ${\rm Pon}(F,E)$ is kind of a too fine way of classifying objects.
For instance, even up to isomorphism the trivial pair over the one-point space $*$  does not have a unique
extension, as $\Gamma(*,{\rm A_{Par}}/{\rm A_{Pon}})={\rm A_{Par}}/{\rm A_{Pon}}\not=\{*\}$.
A more appropriate classification is the classification up to isomorphism in the full
subcategory ${\rm Pon}^{\rm full}(F,E)$ of ${\rm Pon}(B)$. By definition,
over the one-point space there exists only one Pontrjagin duality triple up to isomorphism.

For an ${\rm A_{Par}}$-bundle $P\to B$, let ${\rm Red}^{\rm full}(P)$ be the ``full'' category 
of reductions, i.e. its objects  are reductions to ${\rm A_{Par}}$-bundles 
(\ref{EqThisIsARedNotGreen}) and its morphisms are commutative diagrams
\begin{eqnarray*}
\xymatrix{
P_{r}\ar[rd]\ar[dd]\ar[rr]&&P\ar[dr]\ar[dd]|!{[d];[d]}\hole&\\
&P_{s}\ar[rr]\ar[dd]&&P\ar[dd]\\
B\ar[rd]^{=}\ar[rr]|!{[r];[dr]}\hole^{\ \ \ \ \ =}&&B\ar[dr]^=&\\
&B\ar[rr]^=&&B
}
\end{eqnarray*}
where $P_{r}\to P_{s}$ is an ${\rm A_{Pon}}$-bundle isomorphism over $B$
and $P\to P$ is an ${\rm A_{Par}}$-bundle automorphism over $B$.
There is a direct analogue of Proposition \ref{PropEqOfCatsForRed} 
for the two ``full'' categories:
\begin{prop}
	Let $P\to B$ be an ${\rm A_{Par}}$-bundle, and let 
	$F\to E\to B$ be its associated pair.
	Then there is an equivalence of categories 
	\begin{eqnarray*}
		{\rm Pon}^{\rm full}(F,E)\simeq {\rm Red}^{\rm full}(P).
	\end{eqnarray*}
\end{prop}

We denote by ${\rm Aut}(P)$ the automorphism group of ${\rm A_{Par}}$-bundle 
automorphisms  of $P$ over $B$. The maps $\gamma_{*},\varepsilon_{*}$
of Proposition \ref{PropTheClassOfRedsAndExts} are ${\rm Aut}(P)$-equivariant
with respect to the obvious actions of ${\rm Aut}(P)$. 
Concerning the isomorphism classes of the two ``full'' categories in the proposition above,
we obtain the following result which can be proven by the same means as
Proposition \ref{PropTheClassOfRedsAndExts}.
\begin{prop}\label{PropTheFullExtProblem}
	\begin{itemize}
	\item[i)]
		There are natural bijections 
		$$
		[{\rm Pon}^{\rm full}(F,E)]\cong [{\rm Red}^{\rm full}(P)]
		\cong \frac{{\rm im}(\gamma_{*})}{{\rm Aut}(P)} \cong\frac{{\rm im}(\varepsilon_{*})}{{\rm Aut}(P)}.
		$$
	\item[ii)] If ${\rm A_{Par}}\to {\rm A_{Par}}/{\rm A_{Pon}}$ has local sections, then
		$$
		[{\rm Pon}^{\rm full}(F,E)]\cong [{\rm Red}^{\rm full}(P)]
		\cong \frac{\Gamma(B,P/{\rm A_{Pon}})}{{\rm Aut}(P)}
		\cong \frac{C(P,{\rm A_{Par}}/{\rm A_{Pon}})^{\rm A_{Par}}}{{\rm Aut}(P)}.
		$$
	\end{itemize}
\end{prop}

If $F\to E\to B$ is a pair, it has a 
classifying map $B\to {\rm BA_{\rm Par}}$ into to the classifying 
space of ${\rm A_{\rm Par}}$-principal bundles which is unique 
up to homotopy.
The inclusion (projection) 
${\rm A_{\rm Pon}}\hookrightarrow{\rm A_{\rm Par}},
(a,\phi(a))\mapsto a$
induces a map of the classifying spaces
${\rm BA_{\rm Pon}}\to{\rm BA_{\rm Par}}$,
and an extension of the pair exists if and only if 
there exists a map $B\to{\rm BA_{\rm Pon}}$
such that the diagram
\begin{eqnarray}\label{DiagLiftDingsMittenDrin}
\xymatrix{
&{\rm BA_{\rm Pon}}\ar[d]\\
B\ar[ru]\ar[r]&{\rm BA_{\rm Par}}
}
\end{eqnarray}
commutes up to homotopy.
We have a commutative diagram 
for the set of isomorphism 
classes of triples $[{\rm Pon}(B)]$ over $B$ and 
the set of isomorphism classes 
of pairs $[{\rm Par}(B)]$ over $B$
\begin{eqnarray*}
\xymatrix{
[{\rm Pon}(B)]\ar[d]^{\rm forget}\ar[r]^\cong&[B,{\rm BA_{Pon}}]\ar[d]\\
[{\rm Par}(B)]\ar[r]^\cong&[B,{\rm BA_{Par}}]
},
\end{eqnarray*}
where we put the homotopy classes of maps from $B$ to 
${\rm BA_{Par}}$ or ${\rm BA_{Par}}$ in the right column 
and
the isomorphisms are given by sending an isomorphism class
to the homotopy class of a classifying map.
Let $[F,E]$ denote any element of $[{\rm Par}(B)]$,
then there is a canonical correspondence: 
${\rm forget}^{-1}([F,E])\cong[{\rm Pon}^{\rm full}(F,E)]$,
and in case of $ii)$ of Proposition \ref{PropTheFullExtProblem} 
we end up with a commutative diagram
\begin{eqnarray*}
\xymatrix{
&[{\rm Pon^{full}}(F,E)]\ar@{^(->}[d]\ar[r]^-{\cong}&
\frac{C(P,{\rm A_{Par}}/{\rm A_{Pon}})^{\rm A_{Par}}}{{\rm Aut}(P)}\ar@{^(_->}[d]^{u_{*}}&\\
{\rm forget}^{-1}([F,E])\ar@{->>}[d]\ar@{^(->}[r]\ar@{=}[ru]&
[{\rm Pon}(B)]\ar[d]^{\rm forget}\ar[r]^\cong&[B,{\rm BA_{Pon}}]\ar[d]\\
\{[F,E]\}\ar@{^(->}[r]&[{\rm Par}(B)]\ar[r]^\cong&[B,{\rm BA_{Par}}]
},
\end{eqnarray*}
where $u:P/{\rm A_{Pon}}\to {\rm BA_{Pon}}$ is the classifying 
map of the bundle $P\to P/{\rm A_{Pon}}$ and $u_{*}([\sigma]):= [u\circ\sigma]$.

\begin{exa}
	Let $G:=S^1$ be the circle, so  $\widehat G\cong\Z$.
	Then the inclusion 
	${\rm A_{\rm Pon}}\hookrightarrow{\rm A_{\rm Par}}$
	is a homotopy equivalence.
	It follows that the map ${\rm BA_{\rm Pon}}\to{\rm BA_{\rm Par}}$
	between the classifying spaces is a weak homotopy equivalence.

	So if one considers pairs $F\to E\to B$ over a CW-complex $B$, 
	then up to isomorphism in ${\rm Pon}^{\rm full}(F,E)$ there 
	always exists a unique extension to a Pontrjagin duality triple.
	(Cp. also Example \ref{ExaMitS1Besser})
\end{exa}

Next we try to understand, i.e. simplify, the quotient 
$\frac{C(P,{\rm A_{Par}}/{\rm A_{Pon}})^{\rm A_{Par}}}{{\rm Aut}(P)}$.
Let us assume the pair $F\to E\to B$ has an extension, 
so the ${\rm A_{Par}}$-bundle $P\to B$ has a reduction
to a ${\rm A_{Pon}}$-bundle $Q\to B$.
There are canonical identifications
\begin{eqnarray*}
C(P,{\rm A_{Par}}/{\rm A_{Pon}})^{\rm A_{Par}}\cong 
C(Q,{\rm A_{Par}}/{\rm A_{Pon}})^{\rm A_{Pon}}
\end{eqnarray*}
and
\begin{eqnarray*}
{\rm Aut}(P)\cong C(Q,{\rm A_{Par}}^{\rm ad})^{\rm A_{Pon}},
\end{eqnarray*}
where ${\rm A_{Pon}}\ni b$ acts from the right on ${\rm A_{Par}}^{\rm ad}:={\rm A_{Par}}\ni a$ 
by conjugation: $a\cdot b:=b^{-1} a b$.

Because the normal subgroup 
$(\U(1)\times \widehat G)\subset G\ltimes (\U(1)\times\widehat G)
={\rm A_{Pon}}$
acts trivially on ${\rm A_{Par}}/{\rm A_{Pon}}$,
we have a further identification
$$
C(Q,{\rm A_{Par}}/{\rm A_{Pon}})^{\rm A_{Pon}}\cong C(E,{\rm A_{Par}}/{\rm A_{Par}})^{G},
$$
where we used that $Q/(\U(1)\times \widehat G)\cong E$.
Topologically, the quotient ${\rm A_{Par}}/{\rm A_{Pon}}$
is isomorphic to $C_{*}(G,\U(1))/\widehat G$, where we denoted by 
$C_{*}(G,\U(1))$ the base-point preserving continuous functions.
In fact, we have a diagram of inclusions and quotients of topological spaces
\begin{eqnarray*}
\xymatrix{
G\times\U(1)\ar@{^(->}[r]&G\ltimes(\U(1)\times \widehat G)\ar@{_{(}->}[d]\ar@{->>}[r]&
\ \widehat G\ar@{_{(}->}[d]\ar@<-1ex>@{_(->}[l]\\
G\times\U(1)\ar@{^(->}[r]&G\ltimes C(C,\U(1))\ar@{->>}[r]\ar@{->>}[d]& 
C_{*}(G,\U(1))\ar@<-1ex>@{_(->}[l]\ar@{->>}[d]\\
&\frac{G\ltimes C(C,\U(1))}{G\ltimes (\U(1)\times\widehat G)}\ar[r]^-{\cong}&
\frac{C_{*}(C,\U(1))}{\widehat G}
}.
\end{eqnarray*}
The action of $G$ on ${\rm A_{Par}}/{\rm A_{Pon}}$
induces an action of $G$ on $C_{*}(G,\U(1))/\widehat G$
which turns out to be the action induced by the shift action
of $G$ on $C_{*}(G,\U(1))$:
\begin{eqnarray*}
f\widehat G\cdot g := \big(f(g)^{-1} f(g+\_)\big)\widehat G\in C_{*}(G,\U(1))/\widehat G.
\end{eqnarray*}

Let $Q\ni q\mapsto(g_{q},f_{q})\in {\rm A_{Par}}^{\rm ad}$ be 
an ${\rm A_{Pon}}$-equivariant map, i.e.
\begin{eqnarray*}
(g_{q(h,z,\chi)},f_{q(h,z,\chi)})
&=&(h,z,\chi)^{-1}(g_{q},f_{q})(h,z,\chi)\\
&=&(-h,z^{-1}\langle h,\chi\rangle,-\chi)(g_{q},f_{q})(h,z,\chi)\\
&=&(g_{q},\langle g_{q},\chi\rangle^{-1}f_{q}(h+\_)).
\end{eqnarray*}
In particular, $q\mapsto g_{q}$ is ${\rm A_{Pon}}$-invariant and
$q\mapsto (g_{q},f_{q})$ is $\U(1)$-invariant. 
Furthermore the $\widehat G$-equivariance only contributes by 
the scalar
$q\mapsto \langle g_{q},\chi\rangle^{-1}\in\U(1)$
which itself acts trivially on 
 $C(E,{\rm A_{Par}}/{\rm A_{Pon}})^G$,
 and the mapping $q\mapsto f_{q}(0)^{-1}f_{q}\in C_{*}(G,\U(1))$
 is $G$-equivariant.
So by use of the decomposition 
$G\ltimes C(G,\U(1)) = G\ltimes (\U(1)\times C_{*}(G,\U(1)))$
we can identify 
\begin{eqnarray*}
\frac{C(E,C_{*}(G,\U(1))/\widehat G)^G}
{C(Q,(G\ltimes C(G,\U(1)))^{\rm ad})^{\rm A_{Pon}}}
\cong
\frac{C(E,C_{*}(G,\U(1))/\widehat G)^G}
{C(B,G)\ltimes C(E,C_{*}(G,\U(1)))^G}.
\end{eqnarray*}

\begin{lem} 
	Inside $C(E,C_{*}(G,\U(1))/\widehat G)^G$ 
	the orbits of $C(B,G)\ltimes C(E,C_{*}(G,\U(1)))^G$
	coincide with the orbits of $C(E,C_{*}(G,\U(1)))^G$, 
	so
	\begin{eqnarray*}
		\frac{C(E,C_{*}(G,\U(1))/\widehat G)^G}
		{C(B,G)\ltimes C(E,C_{*}(G,\U(1)))^G}
		\cong
		\frac{C(E,C_{*}(G,\U(1))/\widehat G)^G}
		{C(E,C_{*}(G,\U(1)))^G}.
	\end{eqnarray*}
\end{lem}

\begin{pf}
Let $e\mapsto F_{e}\in C_{*}(G,\U(1))/\widehat G$ be 
$G$-equivariant, and let  $(b\mapsto f_{b})\in C(B,G)$.
For each $F_{e}$ choose (non-continuously) a
$\overline F_{e}\in C_{*}(G,\U(1))$ such that
$\overline F_{e}\widehat G =F_{e}$.
Denote by $p:E\to B$ the projection, then
\begin{eqnarray*}
(f\cdot F)_{e}&=&\big( \overline F_{e}(f_{p(e)})^{-1}\overline F_{e}(f_{p(e)}+\_)\big)\widehat G\\  
&=&\big (\underbrace{ \overline F_{e}(f_{p(e)})^{-1}\overline F_{e}(f_{p(e)}+\_) \overline F_{e}(\_)^{-1}}\big)\widehat G\ F_{e}\\
&&\qquad\qquad\qquad\quad=: d(\overline F_{e})(f_{p(e)},\_),
\end{eqnarray*}
here $d:C(G,\U(1))\to C(G\times G,\U(1))$ is the boundary operator of
group cohomolgy which has kernel $\widehat G$, so it factors (continuously)
\begin{eqnarray*}
\xymatrix{
C(G,\U(1))\ar[d]\ar[r]^-{d}&C(G\times G,\U(1))\\
\frac{C(G,\U(1))}{\widehat G}\ar[ur]_-{\underline d}
},
\end{eqnarray*}
and we have $d (\overline F_{e})=\underline d(F_{e})$.
Therefore $f\cdot F$ differs from $F$ by action of
$$
e\mapsto \underline d(F_{e})(f_{p(e)}+\_)\in C_{*}(G,\U(1))
$$
which is easily checked to be $G$-equivariant. This proves the lemma.
\end{pf}

\begin{cor} 
	If a pair $F\to E\to B$ has an extension, then
	there is a bijection
	\begin{eqnarray*}
		[{\rm Pon}^{\rm full}(F,E)]\cong
		\frac{C(E,C_{*}(G,\U(1))/\widehat G)^G}
		{C(E,C_{*}(G,\U(1)))^G}.
	\end{eqnarray*}
\end{cor}

\begin{exa}\label{ExaMitS1Besser}
	Let $G:=S^1$ be the circle group. Then
	$\widehat G\cong\Z$ and the inclusion  
	$\Z\hookrightarrow C_{*}(S^1,\U(1))$
	is a homotopy equivalence.
	$C_{*}(S^1,\U(1))/\Z$ is $G$-equivariantly
	isomorphic to the null-homotopic functions
	which leads to a $G$-equivariant section
	\begin{eqnarray*}
		C_{*}(S^1,\U(1))/\Z\stackrel{\cong}{\to} 
		C_{*}(S^1,\U(1))_{\rm null}\stackrel{\subset}{\to} C_{*}(S^1,\U(1))
	\end{eqnarray*}
	of the quotient map $C_{*}(S^1,\U(1))\to C_{*}(S^1,\U(1))/\Z$.
	Therefore the quotient 
	$\frac{C(E,C_{*}(S^1,\U(1))/\Z)^{S^1}}
	{C(E,C_{*}(S^1,\U(1)))^{S^1}}$
	consists of a single element only. 
	This means that an extension of a pair
	is unique up to isomorphism.	
\end{exa}

\begin{exa}
Let $G:=\Z$ be the integers, so $\widehat G\cong S^1$ is the circle group.
The quotient map $C_{*}(\Z,\U(1))\to C_{*}(\Z,\U(1))/S^1$
is a canonically trivialisable $S^1$-bundle
\begin{eqnarray}\label{ExaMittendinsontrivialesbuendel}
\xymatrix{
C_{*}(\Z,\U(1)) \ar[d]\ar[r]^-{\cong}&
\prod_{\Z\setminus\{0\}} \U(1)\ar[d]\ar[r]_-{\theta}^-{\cong}& 
{\prod_{\Z\setminus\{0,1\}}\U(1)}\times S^1\ar[d]\\
\frac{C_{*}(\Z,\U(1))}{S^1}\ar[r]^-{\cong}&
\frac{\prod_{\Z\setminus\{0\}} \U(1)}{S^1}\ar[r]^-\cong&
{\prod_{\Z\setminus\{0,1\}} \U(1)}
},
\end{eqnarray}
where $\theta: (z_{k})_{k}\mapsto \left(\left(\frac{z_{k}}{z_{1}^k}\right)_{k}, z_{1}\right)$, and
the action of $S^1\ni t$ on the product $\prod_{\Z\setminus\{0\}}\U(1)\ni (z_{k})_{k}$
is $(z_{k})_{k}\cdot t = (z_{k} t^k)_{k}$.
The action of $\Z\ni n$ on $\prod_{\Z\setminus\{0\}}\U(1)\ni (z_{k})_{k}$
is given by 
$(z_{k})_{k}\cdot n =( \frac{z_{k+n}}{z_{n}})_{k}$. The induced action
of $\Z\ni n$ on 
${\prod_{\Z\setminus\{0,1\}}\U(1)}\times S^1\ni ( (z_{k})_{k},s)$
is
\begin{eqnarray*}
( (z_{k})_{k},s)\cdot n =\left( \left(\frac{z_{k+n} z_{n}^{k-1}}{z_{1+n}^k}\right)_{k},
\frac{z_{1+n}}{z_{n}} s\right).
\end{eqnarray*}
Let $E\to B$ be a $\Z$-bundle, then we are faced with the $\Z$-equivariant lifting
problem: Does there exist a $\psi:E\to S^1$ such that
\begin{eqnarray}\label{Examittendinsonhoheber}
\xymatrix{
&&{\prod_{\Z\setminus\{0,1\}}\U(1)}\times S^1\ar[d]\\
E\ar[rr]^-\varphi\ar@{-->}@/^0.6cm/[rru]^-{\varphi\times\psi}&&{\prod_{\Z\setminus\{0,1\}}\U(1)}
}
\end{eqnarray}
commutes $\Z$-equivariantly for a given $\Z$-equivariant $\varphi$?
The equivariance of $\varphi=(\varphi_{k})_{k}$ gives 
\begin{eqnarray*}
\varphi_{k}(e\cdot n)= \frac{\varphi_{k+n}(e) \varphi_{n}(e)^{k-1}}{\varphi_{1+n}(e)^k},
\qquad e\in E, n\in \Z,
\end{eqnarray*}
 which for $k=2$ implies
\begin{eqnarray}\label{ExaMittendrinsonphizwo}
\varphi_{2+n}(e)=\varphi_{2}(e\cdot n)\varphi_{1+n}(e)^2\varphi_{n}(e)^{-1}.
\end{eqnarray}
As $\varphi_{0}=\varphi_{1}=1$, this shows that all $\varphi_{k}$
can be expressed in terms of $\varphi_{2}$, and, conversely, 
any $\varphi_{2}:E\to \U(1)$ determines a $\varphi$ via 
(\ref{ExaMittendrinsonphizwo}). 
If $\psi$ exists, then by equivariance
$\psi(e\cdot n)=\varphi_{1+n}(e)/\varphi_{n}(e)\psi(e)$. In particular
\begin{eqnarray}\label{ExaMittendrinsonpsi}
\psi(e\cdot 1)=\varphi_{2}(e)\psi(e)
\end{eqnarray}
(note again $\varphi_{0}=\varphi_{1}=1$).
The solvability of the lifting problem  (\ref{Examittendinsonhoheber}) 
is equivalent to to solvability of (\ref{ExaMittendrinsonpsi}),
and it depends on $E$ whether or not this equality can be always solved:
\begin{itemize}

\item[i)] Let $E:=\R\to T=:B$ be the universal covering of the 1-torus $T:=\R/\Z$.
Then (\ref{ExaMittendrinsonpsi}) can always be solved by the formula
\begin{eqnarray*}
\psi(x):=
\begin{cases}
\prod_{k=0}^{n-1}\varphi_{2}(\varepsilon+k)\eta(\varepsilon),\text{ if } x=\varepsilon+n,\varepsilon\in[0,1), n\in\{0,1,2,\dots\}\\
\prod_{k=1}^{n}\varphi_{2}(\varepsilon-k)^{-1}\eta(\varepsilon),\text{ if } x=\varepsilon-n,\varepsilon\in[0,1),n\in\{1,2,\dots\}
\end{cases},
\end{eqnarray*}
where $\eta:[0,1]\to\U(1)$ is any continuous map satisfying $\eta(1)=\varphi_{2}(0)\eta(0)$.
Therefore the pair
\begin{eqnarray*}
\R\times\U(1)\to \R\to T
\end{eqnarray*}
has a unique extension.

\item[ii)] Let $E:= T\times R\to T^2=:B$.
Let $\varphi_{2}:T\times \R\to T\stackrel{k}{\to}\U(1)$ be 
a map with winding number $k\in\Z$.
The equality $\psi(t,x+1)=\varphi_{2}(t,x)\psi(t,x)$
can be solved iff $k=0$, and each $k\in\Z$ defines a 
class in the quotient 
$\frac{C(E,C_{*}(G,\U(1))/\widehat G)^G}
		{C(E,C_{*}(G,\U(1)))^G}$. 
Therefore the pair 
\begin{eqnarray*}
T\times\R\times\U(1)\to T\times\R\to T^2
\end{eqnarray*}
has $\Z$-many extensions.
\end{itemize}

\end{exa}

\begin{exa}
Let $G:={\rm O}(1)\cong \widehat G$ be the group with two elements.
We have a commutative diagram
\begin{eqnarray*}
	\xymatrix{
	C_{*}({\rm O}(1),\U(1))\ar[r]^-{\cong}\ar[d]& \U(1)\ar[d]^2\\
	\frac{C_{*}({\rm O}(1),\U(1))}{\widehat{{\rm O}(1)}}\ar[r]^-{\cong}&\U(1)
	},
\end{eqnarray*}
where $2:\U(1)\to \U(1)$ is the two-fold covering.
The induced action by the non-trivial element of $G$
on the two occurring $\U(1)$ is just complex conjugation 
(inversion).

Let $E:=S^n\to \R{\rm P}^n=:B$ be the two-fold covering,
so the non-trivial element of ${\rm O}(1)$ acts by
mapping a point of $S^{n}$ to its antipodal point.

Then the two constant maps $S^n\to \U(1): x\mapsto \pm 1$ 
to the fixpoints of the conjugation action are 
$G$-equivariant. Clearly, $x\mapsto +1$ can be lifted $G$-equivariantly:
\begin{eqnarray*}
	\xymatrix{
	&& \U(1)\ar[d]^2\\
	S^{n}\ar[rr]^-{x\mapsto +1}\ar[urr]^-{x\mapsto+1}&&\U(1)
	}
\end{eqnarray*}
which is not the case for $x\mapsto -1$. 
In fact, these two maps represent two classes such that
\begin{eqnarray*}
\frac{C(S^n, C_{*}({\rm O}(1),\U(1))/{\widehat{\rm O}(1))^{{\rm O}(1)}}}
{C(S^n,C_{*}({\rm O}(1),\U(1)))^{{\rm O}(1)}}
\cong\{ [x\mapsto +1],[x\mapsto-1]\}
\end{eqnarray*}
Therefore  the pair
\begin{eqnarray*}
S^n\times\U(1)\to S^n\to\R{\rm P}^n
\end{eqnarray*}
has two non-isomorphic extensions
(cp. Example \ref{ExaNTRPWUTP}).

\end{exa}

Let us consider the case of ring pairs $F_{0}\to E_{0}\to B$, so
$E_{0}=B\times G$ is the trivial bundle.
In this case $G$-equivariant maps on $E_{0}$
can be identified with (non-equivariant) maps 
on the base, so 
\begin{eqnarray}\label{EqAufeinmaleinquader}
[{\rm Pon}^{\rm full}(F_{0},E_{0})]\cong 
\frac{C(B,C_{*}(G,\U(1))/\widehat G)}{C(B,C_{*}(G,\U(1)))}.
\end{eqnarray}

\begin{exa}
Let $G:=\Z$ be the integers, so $\widehat G\cong \Z$. 
As (\ref{ExaMittendinsontrivialesbuendel})
is a trivial bundle the quotient (\ref{EqAufeinmaleinquader}) 
is trivial. Therefore all ring pairs 
$F_{0}\to B\times \Z\to B$ admit a unique extensions.
\end{exa}

Let us introduce some self-explaining terminology.
\begin{defi}
	A {\bf dual} of a pair $F\to E\to B$ is a dual pair 
	$\widehat F\to\widehat E\to B$ such that the diagram 
	\begin{eqnarray}\label{DiatheExtPofPaDP}
	\resizebox{!}{0.2cm}{
		\xymatrix{
		F\ar[dr]&&&&\widehat F\ar[dl]\\
		&E\ar[dr]&&\widehat E\ar[dl]\\
		&&B&&
		}}
	\end{eqnarray}
	admits an extension. In that case we also call the 
	pair a {\bf dual} of the dual pair. 
\end{defi}

For a pair and a dual of it the discussion of the extension problem
(\ref{DiatheExtPofPaDP}) is to answer the question of how many 
different ways there are for being dual to each other.

Let ${\rm Pon}^{\rm top}(F,E,\widehat F,\widehat E)$ be the full subcategory
of ${\rm Pon}(E,\widehat E)$ (p. \pageref{PageOfPonTrip}) whose
objects are extensions of (\ref{DiatheExtPofPaDP}). 
By Corollary \ref{CorTheClassOfPon}, we know that the isomorphism
classes of ${\rm Pon}(E,\widehat E)$ are a 
$\check H^1(B,\underline{\U(1)})$-torsor.
The projections $p:E\to B$ and $\hat p:\widehat E\to B$ induce
homomorphisms
$p^*:\check H^1(B,\underline{\U(1)})\to  \check H^1(E,\underline{\U(1)})$
and $\hat p^*:\check H^1(B,\underline{\U(1)})\to \check H^1(\widehat E,\underline{\U(1)})$.
The subgroup $N(E,\widehat E):={\rm ker}(p^*)\cap {\rm ker}(\hat p^*)
\subset \check H^1(B,\underline{\U(1)})$
still acts on the set 
$$[{\rm Pon}^{\rm top}(F,E,\widehat F,\widehat E)]\subset [{\rm Pon}(E,\widehat E)]
$$ of isomorphism
classes of ${\rm Pon}^{\rm top}(F,E,\widehat F,\widehat E)$,
and this action is still free and transitive:

\begin{cor}
	The set of isomorphism classes $[{\rm Pon}^{\rm top}(F,E,\widehat F,\widehat E)]$
	is a $N(E,\widehat E)$-torsor.
\end{cor}

Let us denote by ${\rm Pon^{full}}(F,E,\widehat F,\widehat E)$
the full subcategory of ${\rm Pon}(B)$ whose objects are the extensions 
of (\ref{DiatheExtPofPaDP}).
The inclusion functor
${\rm Pon}^{\rm top}(F,E,\widehat F,\widehat E)\hookrightarrow
{\rm Pon^{full}}(F,E,\widehat F,\widehat E)$ induces
a surjection on isomorphism classes
$$
[{\rm Pon}^{\rm top}(F,E,\widehat F,\widehat E)]\twoheadrightarrow
[{\rm Pon^{full}}(F,E,\widehat F,\widehat E)]
$$
which is $N(E,\widehat E)$-equivariant for the 
descended action of $N(E,\widehat E)$ on  
$[{\rm Pon^{full}}(F,E,\widehat F,\widehat E)]$.
In fact, this action is well-defined and therefore still transitive, 
and we shall describe its stabiliser group next.

The cup product 
$\cup:\check H^0(B,\underline { G})\times \check H^1(B,\underline{\widehat G})
\to\check H^1(B,\underline{\U(1)})$ and
the class $[\widehat E]\in \check H^1(B,\underline {\widehat G})$ 
of $\widehat E\to B$ define a map
$\_\cup[\widehat E] :\check H^0(B,\underline G)\to\check H^1(B,\underline{\U(1)})$.
Its image ${\rm im}(\_\cup[\widehat E] )$ is contained in ${\rm ker}(\hat p^*)$,
as the pullback $\hat p^*\widehat E\to \widehat E$ is trivialisable. 
Dually, ${\rm im}( [E]\cup \_ )\subset {\rm ker}(p^*)$ for 
the cup product 
$\cup:\check H^1(B,\underline { G})\times \check H^0(B,\underline{\widehat G})
\to\check H^1(B,\underline{\U(1)})$.
Let us denote by 
$M(E,\widehat E)$ the following subgroup
\begin{eqnarray*}
M(E,\widehat E)
&:=&\big({\rm im}(\_\cup[\widehat E] )\cap {\rm ker}(p^*)\big)+
\big({\rm im}([ E]\cup\_ )\cap {\rm ker}(\hat p^*)\big)\\
&=&\big( {\rm im}(\_\cup[\widehat E] )+ {\rm im}([ E]\cup\_ )\big)\cap N(E,\widehat E)\\
&\subset&N(E,\widehat E).
\end{eqnarray*}

\begin{prop}\label{Propnendickenstabilisator}
	The stabiliser of the $N(E,\widehat E)$-action on  
	$[{\rm Pon^{full}}(F,E,\widehat F,\widehat E)]$
	is $M(E,\widehat E)$, so
	$[{\rm Pon^{full}}(F,E,\widehat F,\widehat E)]$
	is a $N(E,\widehat E)/M(E,\widehat E)$-torsor.
\end{prop}

\begin{pf}
Note first that we have canonical identifications
\begin{eqnarray*}
\check H^0(B,\underline{ G})=C(B, G)\cong {\rm Aut}(E),
\text{ and }
\check H^0(B,\underline{\widehat G})=C(B,\widehat G)\cong {\rm Aut}(\widehat E),
\end{eqnarray*}
where ${\rm Aut}(E), {\rm Aut}(\widehat E)$ are the bundle 
automorphisms over $B$. Pick an arbitrary $\varphi\in C(B,G)$,
we denote by $\varphi_{E}$ the corresponding bundle automorphism
of $E$. Dually,
$\hat \varphi_{\widehat E}$ is the bundle automorphism
of $\widehat E$ corresponding to some $\hat \varphi\in C(B,\widehat G)$.
An explicit bundle $\varphi\cup\widehat E$
representing the class $\varphi\cup [\widehat E]$ is
given by the pullback
\begin{eqnarray}\label{DiagTheCupAsIcecream}
\xymatrix{
\varphi\cup \widehat E\ar[r]\ar[d] &F_{\widehat E}\ar[d]\\
B\ar[r]^-{{\rm id}\times\varphi}& B\times G
},
\end{eqnarray}
where $F_{\widehat E}:= \widehat E\times_{\widehat G}(G\times\U(1))$ as before.
We have a canonical isomorphism $(-\varphi)\cup \widehat  E^{\rm op}\cong \varphi\cup\widehat E$.
Dually, we define the pullback $E\cup\hat\varphi$ which represents 
the class $[E]\cup\hat\varphi$.
Now, let 
\begin{eqnarray}\label{DiagTripIndubySattelites1}
	\resizebox*{!}{2,9cm}{$
	\begin{array}{c}
	\xymatrix{
	&F\times_B \widehat E\ar[rd]\ar[ld]\ar[rr]^\kappa&
	&E\times_B \widehat F\ar[rd]\ar[ld]&\\
	F\ar[rd]&&E\times_B \widehat E\ar[rd]\ar[ld]&&\widehat F\ar[ld]\\
	&E\ar[rd]^p&&\widehat E\ar[ld]_{\hat p}&\\
	&&B&&
	}
	\end{array}
	$}
\end{eqnarray}
be a Pontrjagin duality triple, and consider the induced triple
given by tensoring the $\U(1)$-bundles with the $\U(1)$-bundle
$(\varphi\cup \widehat E)\otimes (E\cup\hat \varphi)\to B$:
\begin{eqnarray}\label{DiagTripIndubySattelites2}
	\resizebox*{14cm}{3,7cm}{$
	\begin{array}{c}
	\xymatrix@!C{
	&F\otimes p^*(\varphi\cup \widehat E\otimes  E\cup\hat \varphi)\times_B \widehat E\ar[rd]\ar[ld]\ar[rr]^{\kappa\otimes {\rm flip}}&
	&E\times_B \widehat F\otimes \hat p^*(\varphi \cup\widehat E\otimes E\cup\hat\varphi)\ar[rd]\ar[ld]&\\
	F\otimes p^*(\varphi\cup \widehat E\otimes  E\cup\hat \varphi)\ar[rd]&
	&E\times_B \widehat E\ar[rd]\ar[ld]&
	&\widehat F\otimes \hat p^*(\varphi\cup\widehat E\otimes  E\cup\hat \varphi)\ar[ld]\\
	&E\ar[rd]^p&&\widehat E\ar[ld]_{\hat p}&\\
	&&B&&
	}
	\end{array}
	$}.
\end{eqnarray}
To prove the proposition it is sufficient to prove that
$(-\varphi)_{E}:E\to E$ and $(-\hat\varphi)_{\widehat E}:\widehat E\to \widehat E$ can be extended 
to an isomorphism from
(\ref{DiagTripIndubySattelites1}) to (\ref{DiagTripIndubySattelites2})
in ${\rm Pon}(B)$.

First note that we have a canonical trivialisation
$\hat p^*(\varphi\cup\widehat E)\cong \widehat E\times\U(1)$
which is induced by the canonical map
$$
\xymatrix{
&& F_{\widehat E}\ar[d]\\
\widehat E\ar@/^0.5cm/[rru]^-{[{\rm id}\times(\varphi\circ p)\times 1]}\ar[r]^{\hat p}&B\ar[r]^-{{\rm id}\times\varphi}&B\times G
}
$$
Dually, we have a canonical identification $p^*(E\cup\hat\varphi)\cong E\times\U(1)$.
We define a bundle isomorphism over $(-\varphi)_{E}$ by
$$
\xymatrix{
F\otimes p^*(\varphi\cup\widehat E)\ar[r]^-{\cong}\ar[dd]
&F\otimes p^*((-\varphi)\cup\widehat E^{\rm op})\ar[r]\ar@{.>}[d]& F\ar[dd]\\
&F\times_{B}F_{\widehat E^{\rm op}}\ar[ur]_{\varrho_{\kappa}}\\
E\ar[rr]^{(-\varphi)_{E}}&&E
},
$$
where $\varrho_{\kappa}$ is from (\ref{EqWerhettgeadcachtdat}),
and the dottet arrow sends an element
$[x,(e,[\hat e^{\rm op},-\varphi(p(e)),t)]]$ to 
$(x,[\hat e^{\rm op},-\varphi(p(e)),t)] )$ which is not 
well-defined, but its composition with $\varrho_{\kappa}$ is.
Dually, we obtain a bundle isomorphism 
$\widehat F\otimes \hat p^*(E\cup\hat\varphi)\to \widehat F$
over 
$(-\hat\varphi)_{\widehat E}$. 
It is then a straight forward calculation to show that
the diagram
$$
\xymatrix{
F\otimes p^*(\varphi\cup E\otimes  E\cup\hat \varphi)\times_B \widehat E
\ar[d]^\cong\ar[rr]^{\kappa\otimes {\rm flip}}
&&
E\times_B \widehat F\otimes \hat p^*(E\cup\hat\varphi\otimes E\cup\hat\varphi)\ar[d]^\cong \\
F\times_{B}\widehat E\ar[rr]^\kappa&&E\times_{B}\widehat F
}
$$
commutes. This proves the proposition.
\end{pf}

\begin{rem}
Let ${\rm Pon}^{\rm full}(E,\widehat E)$ be the full subcategory
of ${\rm Pon}(B)$ consisting of extensions of 
$E\to B\leftarrow \widehat E$.
In fact, the proof of Proposition \ref{Propnendickenstabilisator} 
shows that 
$
[{\rm Pon}^{\rm full}(E,\widehat E)]
$
is a torsor for the group $\check H^1(B,\underline{\U(1)})/({\rm im}(\_\cup [\widehat E])+{\rm im}([E]\cup\_))$.
The set $[{\rm Pon}^{\rm full}(E,\widehat E)]$
decomposes into a disjoint union
\begin{eqnarray*}
[{\rm Pon}^{\rm full}(E,\widehat E)]=\coprod [{\rm Pon}^{\rm full}(F,E,\widehat F,\widehat E)]
\end{eqnarray*}
such that each component is a 
torsor for the group
$
N(E,\widehat E)/M(E,\widehat E).
$
\end{rem}

\begin{exa}
Let $B:=\Sigma$ be a connected, two-dimensional, orientable, closed manifold,
then $\check H^2(\Sigma,\underline{\U(1)})\cong H^3(\Sigma,\Z)=0$,
so all choices of bundles $E\to \Sigma\leftarrow\widehat E$
can be extended to a Pontrjagin dualitiy triple.
Let $G:=S^1$ be the circle group, so $\widehat G\cong \Z$,
and let
$E_{1}\to \Sigma$ be an $S^1$-bundle 
representing a generator of 
$\check H^1(\Sigma,\underline{\U(1)})\cong H^2(\Sigma,\Z)\cong\Z$.
Then 
$[E_{1}]\cup\_ :\check H^0(\Sigma,\underline\Z)\to \check H^1(\Sigma,\U(1))$
is surjective. In fact, it is the identity map after identifying 
$\check H^0(\Sigma,\underline\Z)\cong \Z$ and 
$\check H^{1}(\Sigma,\underline{\U(1)})\cong \Z$.
Therefore 
$$\check H^1(\Sigma,\underline{\U(1)})/
({\rm im}(\_\cup [\widehat E])+{\rm im}([E_{1}]\cup\_))=0$$
for any $\Z$-bundle $\widehat E\to\Sigma$,
and all diagrams $E_{1}\to\Sigma\leftarrow \widehat E$
have (up to isomorphism) a unique extension in ${\rm Pon}^{\rm full}(E_{1},\widehat E)$.
\end{exa}

\begin{exa}
Let $T:=\R/\Z$, and let $B:=T^2$ be the two-dimensional torus.
Let $G:=S^1$ be the circle group,
so $\widehat G\cong\Z$.
A $\U(1)$-bundle representing the generator of $H^2(T^2,\Z)$
is 
\begin{eqnarray*}
\xymatrix{
E_{1}:=F_{\R}:=\R\times_{\Z}(S^1\times\U(1))\ar[d]\\
T^2=T\times S^1
}.
\end{eqnarray*}
Let $\widehat E_{m,n}\to T^2$ be the 
$\Z$-bundle given by pullback
\begin{eqnarray*}
\xymatrix{
\widehat E_{m,n}\ar[r]\ar[d]& \R\ar[d]\\
T^2\ar[r]^{f_{m,n}}&T
}
\end{eqnarray*}
along the map $f_{m,n}: (x,y)\mapsto mx+ny$
which represents 
$(m,n)\in\Z\oplus\Z\cong H^1(T^2,\Z)$.
For an integer $k$ consider the 
extension problem 
\begin{eqnarray}\label{DiagMyFavExtMineralwasser}
\xymatrix{
E_{1}^{\otimes k}\ar[dr]&& \widehat E_{m,n}\ar[dl]\\
&T^2&
}.
\end{eqnarray}
The image of $[E_{1}^{\otimes k}]\cup\_ $ in
$\check H^{1}(T^2,\underline{\U(1)})$ is 
the subgroup corresponding to 
$k\Z\subset \Z\cong  H^{2}(T^2,\Z)$.
To understand the image of $\_\cup[\widehat E_{m,n}]$
in $\check H^{1}(T^2,\underline{\U(1)})$ we have to 
consider the diagram of pullbacks
according to (\ref{DiagTheCupAsIcecream})
\begin{eqnarray*}
\xymatrix{
\varphi\cup \widehat E_{m,n}\ar[d]\ar[r]& \widehat E_{m,n}\times_{\Z}(S^1\times\U(1))
\ar[d]\ar[r]&
\R\times_{\Z}(S^1\times\U(1))\ar[d]\\
T^2\ar[r]^{{\rm id}\times\varphi}\ar@/_1cm/[rr]^{f_{m,n}\times\varphi}& T^2\times S^1\ar[r]^{f_{m,n}\times{\rm id}}&
T\times S^1
}.
\end{eqnarray*}
For $\varphi={\rm pr_{2}}: (x,y)\mapsto y$ 
it follows that ${\rm pr_{2}}\cup \widehat E_{m,n}$ is isomorphic to 
$E_{1}^{\otimes m}$, and  for $\varphi ={\rm pr_{1}}:(x,y)\to x$ it follows that
${\rm pr_{1}}\cup \widehat E_{m,n}$ is isomorphic to $E_{1}^{\otimes (-n)}$.
Therefore, if any two of the three numbers $k, m$ and $n$ are coprime, then
$$
{\rm im}([E_{1}^{\otimes k}]\cup\_ )+
{\rm im}(\_\cup[\widehat E_{m,n}])= \check H^1(T^2,\underline{\U(1)})
$$
which means that in this case the extension problem
(\ref{DiagMyFavExtMineralwasser}) has 
a solution which is unique (up to isomorphism in 
${\rm Pon}^{\rm full}(E_{1}^{\otimes k},\widehat E_{m,n})$).
\end{exa}

\section{The Fourier transform}
\label{SecTheFourierTransform}

In this section we show that the Fourier transform 
implements a functor from Pontrjagin duality triples
to the category of tuples consisting of two isomorphic Hilbert 
$C^*$-modules. 
Let us explain what is meant by this.

Let ${\rm Pon}(B)$ denote the category
of Pontrjagin duality triples over a space
$B$, and let ${\rm Pon}$ denote the 
associated total category of all Pontrjagin 
duality triples over all spaces. I.e.
a morphism from a Pontrjagin
duality triple $X$ over $B$ to 
a Pontrjagin duality triple $X'$ over
$B'$ consists of a continuous map
$f:B\to B'$ and an isomorphism 
$X\to f^*X'$ in ${\rm Pon}(B)$.
Let $(0\cong1)$ denote the category 
with two objects $0,1$ and two non-trivial 
morphisms $0\to 1$ and $1\to 0$ which are
inverse to each other.
Let ${\rm Hilb}\text{-}C^*\text{-}{\rm mod}$ denote
the category of Hilbert $C^*$-modules. I.e. a
morphism from  $(H,C)$ to $(H',C')$
consists of a continuous linear map $H\to H'$ of 
Banach spaces  and a morphism of $C^*$-algebras 
$C\to C'$ such that the obvious conditions are satisfied.
In this section we define  a bifunctor
\begin{eqnarray*}
({\rm Pon})^{\rm op}\times (0\cong1)\to {\rm Hilb}\text{-}C^*\text{-}{\rm mod},
\end{eqnarray*}
equivalently,  a functor   
\begin{eqnarray*}
({\rm Pon})^{\rm op}\to ({\rm Hilb}\text{-}C^*\text{-}{\rm mod})^{(0\cong1)},
\end{eqnarray*}
where the target category is a functor category.
This functor will be constructed such that the trivial Pontrjagin duality triple over
the one-point space
is mapped to the classical isomorphism of 
Pontrjagin duality:
\begin{eqnarray*}
	\boldsymbol{
	\resizebox*{5cm}{3cm}{$
	\left(
	\begin{array}{c}	
		\xymatrix{
		&G\times \widehat G\times\U(1)\ar[rd]\ar[ld]\ar[rr]^{\pi}&
		&G\times \widehat G\times\U(1)\ar[rd]\ar[ld]&\\
		G\times \U(1)\ar[rd]&&G\times \widehat G\ar[rd]\ar[ld]&&
		\widehat G\times\U(1)\ar[ld]\\
		&G\ar[rd]&&\widehat G\ar[ld]&\\
		&&\ast&&
	}\end{array}
	\right)$
	}}
	\mapsto
	\Big(\big( C^*(G),C^*(G) \big)\cong \big(C_0(\widehat G),C_0(\widehat G)\big)\Big).
\end{eqnarray*}
Here the group $C^*$-algebra of $G$ and the at infinity vanishing 
functions on $\widehat G$ are understood as Hilbert $C^*$-modules
over themself.

Let us now introduce the Fourier transform based 
on a  Pontrjagin duality triple. Let 
\begin{eqnarray}\label{DiagZehntausendPDT}
	\xymatrix{
	&F\times_B \widehat E\ar[rd]\ar[ld]\ar[rr]^\kappa&
	&E\times_B \widehat F\ar[rd]\ar[ld]&\\
	F\ar[rd]&&E\times_B \widehat E\ar[rd]\ar[ld]&&\widehat F\ar[ld]\\
	&E\ar[rd]&&\widehat E\ar[ld]&\\
	&&B&&
	}
\end{eqnarray}
be any Pontrjagin duality triple. We denote by 
$F^\C$ and $\widehat F^\C$ the associated line
bundles of $F$ and $\widehat F$, respectively.
We define the {\bf Fourier
transform based on diagram (\ref{DiagZehntausendPDT})}
to be a map
$$ 
{\hat{}}:\Gamma_c(E,\widehat F^\C)\to \Gamma_0(\widehat E,\widehat F^\C).
$$
Here $\Gamma_c(E,\widehat F^\C)$ and $\Gamma_0(\widehat E,\widehat F^\C)$
are the bounded continuous section which have fibre-wise compact support 
or which vanish at fibre-wise infinity as indicated by $\Gamma_c$
and $\Gamma_0$, respectively.
The top isomorphism $\kappa$ induces an 
isomorphism of line bundles
 $\kappa^\C: F^\C\times_B\widehat E\to E\times_B \widehat F^\C$. 
For a fibre-wise compactly supported section $\gamma: E\to F^\C$ 
we define its Fourier transform $\hat\gamma:\widehat E\to \widehat F^\C$ by 
\begin{eqnarray}\label{EqTheFTofSections}
	\hat \gamma(\hat e):=\int_G  {\rm pr}_{\widehat F^\C}\big(
	\kappa^\C(\gamma(e\cdot h),\hat e)\big)\ dh,
\end{eqnarray}
where $e\in E$ is any point over the image of $\hat e$ in $B$.
Note that $\hat\gamma$ is a well-defined  section 
as ${\rm pr}_{\widehat F^\C}(
\kappa^\C(\gamma(e\cdot h),\hat e)$ stays 
in the same fibre $\widehat F^\C|_{\hat e}$
independent of $h\in G$ .

As an example take the Pontrjagin duality triple 
 (\ref{DiagTheFirstExaOfAPDT})
given by a $\widehat G$-bundle $\widehat E\to B$:
\begin{eqnarray}\label{DiagBayrischerWald}
\xymatrix{
&F_{\widehat E^{\rm op}}\times_B \widehat E\ar[rd]\ar[ld]\ar[rr]^{\kappa_{\widehat E}}&
&\ar[rd]\ar[ld](B\times G)&\hspace{-1,6cm} \times_B \widehat E\times \U(1)\\
F_{\widehat E^{\rm op}}\ar[rd]&&(B\times G)\times_B \widehat E\ar[rd]\ar[ld]&&\widehat E\times \U(1) \ar[ld]\\
&B\times G\ar[rd]&&\widehat E\ar[ld]&\\
&&B&&
}.
\end{eqnarray}
The Fourier transform based on (\ref{DiagBayrischerWald})
maps a section $\alpha :B\times G\to F^\C_{\widehat E^{\rm op}}$  
to a section $\hat \alpha:\widehat E\to \widehat E\times \C$
of the trivial line bundle,
so it is a map 
$$
\Gamma_{c}(B\times G,F^\C_{\widehat E^{\rm op}})\to C_0(\widehat E),
$$
where $C_0(\widehat E)$ are the bounded continuous functions 
vanishing at fibre-wise infinity.
\\

The intention of this section is to embed the Fourier transform 
based on a Pontrjagin duality triple into a $C^*$-algebraic context.

We start with the definition of the $C^*$-algebra $C^*(E_0,F_0)$
of a ring pair  $F_0\to E_0\to B$, $F_0=F_{\widehat E}, E_0= B\times G$.
As we observed in section \ref{SecRingPairs}, 
a ring pair can be thought of as a $\U(1)$-central extension of 
groupoids 
\begin{eqnarray*}
\xymatrix{
B\!\times\!\U(1)\ \ar@<2pt>[d]\ar@<-2pt>[d]\ar@{>->}[r]&F_0\ar@{->>}[r]\ar@<2pt>[d]\ar@<-2pt>[d] 
&E_0\ar@<2pt>[d]\ar@<-2pt>[d]\\
B\ar[r]^=&B\ar[r]^=&B
}.
\end{eqnarray*}
In short,  the $C^*$-algebra $C^*(E_0,F_0)$ of the ring pair  
is  the $C^*$ algebra of this central extension \cite[Section 3]{TXL}. 

Let us spell out what this means. 
Denote by $F_0^\C:=F_0\times_{\U(1)}\C$ the associated line bundle.
The multiplication map $\mu$ of the ring pair induces 
an associated map $\mu^\C:F^\C_0\times_{B} F^\C_0\to F^\C_0$
by $\mu^\C([x,z],[x',z']):=[\mu(x,x'),z z']$, and we
let $[x,z]^\ddag:=[{x^\dag},\overline{z}]$, where
$ x^\dag\in F_0$ is as in Remark \ref{RemRingPairIsGroupBndl}
and $\overline z$ is the complex conjugate of $z\in\C$.
Let $\Gamma_{c}(E_0,F_0^\C)$ denote the vector space 
of bounded continuous sections with compact support along
the fibres.
We define a convolution for two such sections $\alpha,\beta$ by
\begin{eqnarray*}
(\alpha \ast_\mu \beta)(b,g):=\int_G \mu^\C(\alpha(b,g-h),\beta(b,h))\ dh, 
\end{eqnarray*}
a star operation by
\begin{eqnarray*}
\alpha^\ast(b,g):= (\alpha(b,-g))^\ddag,
\end{eqnarray*}
and a norm by taking the supremum over all fiber-wise
$L^1$-norms, i.e.
\begin{eqnarray*}
||\alpha||_{(\infty,1)}:=\sup_{b\in B} \int_G |\alpha(b,g)|\ dg,
\end{eqnarray*}
where $|[x,z]|:=|z|$. 
The completion of $\Gamma_{c}(E_0,F_0^\C)$ with respect 
to this norm yields a Banach $^*$-algebra. We call its enveloping 
$C^*$-algebra $C^*(E_0,F_0)$  the $C^*$-algebra of 
the ring pair $F_0\to E_0\to B$. 
\\

Next we give the definition of the Hilbert $C^*$-module
$H(E,F)$ of a module pair  $F\to E\to B$  over $F_0\to E_0\to B$.
We observed in section \ref{SecModulePairs}
that we can interpret a module pair as a $\U(1)$-equivariant
Morita self-equivalence 
\begin{eqnarray*}
	\Big(
	\xymatrix{
	F_0 &
	\ar@{}[l]|{\ \ \resizebox{!}{0,6cm}{\rotatebox {90}{{\rotatebox {90}{
	{\rotatebox {90}{$\circlearrowright$}}}}}}}
	F\ar[dl]\ar[dr] \ar@{}[r]|{\resizebox{!}{0,6cm}{\rotatebox {90}{$\circlearrowleft$}}}
	& F_0\\
	B\ar@{<-}@<2pt>[u]\ar@{<-}@<-2pt>[u]&&B\ar@{<-}@<2pt>[u]\ar@{<-}@<-2pt>[u]
	}\Big)\ {\rotatebox {90}{$\circlearrowleft$}}\ \U(1).
\end{eqnarray*}
of the groupoid central extension given by the ring pair.
The Hilbert $C^*$-module $H(E,F)$ is the Hilbert $C^*$-module of 
this Morita equivalence. It is a module over $C^*(E_0,F_0)$.

Again we spell out some of this construction.
Denote by $\varrho:F\times_B F_0\to F$ the action map, and
denote by $F^\C$ the associated line bundle to $F$.
We use the same formula we have used to define $\mu^\C$ 
to define $\varrho^\C:F^\C\times_B F_0^\C\to F^\C$.
The vector space of bounded, fibre-wise compactly supported  
sections $\Gamma_{c}(E,F^\C)$ 
has  the structure of a  $\Gamma_{c}(E_0,F_0^\C)$-right module
by
\begin{eqnarray*}
	(\gamma \ast_\varrho \alpha)(e):=
	\int_G \varrho^\C(\gamma(e\cdot(-h)),\alpha(b,h))\ dh.
\end{eqnarray*}
Here $b\in B$ is the image of $e\in E$ by $E\to B$.
We define  $\sigma_\varrho^\C:F^\C\times_B F^\C\to F_0^\C$
by $\sigma_\varrho^\C([y,z],[y',z']):=[\sigma_\varrho(y,y'),\overline z z']$,
then $\sigma_\varrho^\C([y',z'],[y,z])=\sigma_\varrho^\C([y,z],[y,'z'])^\ddag$
(cp. Remark \ref{RemAboutModPairs}).
For $\gamma,\delta\in \Gamma_{c}(E,F^\C)$ we have 
a $\Gamma_{c}(E_0,F_0^\C)$-valued inner product by
\begin{eqnarray*}
\langle \gamma,\delta\rangle_{c}(b,g):=
\int_G \sigma_\varrho^\C( \gamma(e\cdot h),\delta(e\cdot (h+g)))\ dh \in F_0^\C|_{b\times g}.
\end{eqnarray*}
One can check that this is a sesquilinear form which is 
$\Gamma_{c}(E_0,F_0^\C)$-linear, i.e.
$
\langle \gamma,\delta\rangle_{c}\ast_\mu \alpha= \langle\gamma,\delta\ast_\varrho\alpha\rangle_c,
$
which  is anti-symmetric, i.e. 
$\langle\delta,\gamma\rangle_{c}=\langle\gamma,\delta\rangle_{c}^\ast$,
which is positive, i.e. $\langle \gamma,\gamma\rangle_{c} \ge 0\in C^*(E_0,F_0^\C)$,
and which is definite, i.e. $\langle\gamma,\gamma\rangle_{c}=0$  implies $\gamma=0$.
We define a norm on $\Gamma_{c}(E,F^\C)$ by
\begin{eqnarray*}
	||\gamma||_{H(E,F)} :=
	\Big( ||\langle\gamma,\gamma\rangle_{c}||_{C^*(E_0,F_0)}\Big)^{\frac{1}{2}},
\end{eqnarray*}
and $H(E,F)$ is the completion 
of $\Gamma_{c}(E,F^\C)$ with respect to this norm.

\begin{lem}
	The maps $\ast_\varrho$ and $\langle\  ,\ \rangle_c$ as defined above 
	extend together with its properties
	to $H(E,F)$ and $C^*(E_0,F_0)$ such that
	$H(E,F)$ becomes a full Hilbert $C^*$-module over $C^*(E_0,F_0)$. 
\end{lem}

\begin{pf}
	The proof is standard.
\end{pf}

Now consider a Pontrjagin duality triple
	\begin{eqnarray*}
	\resizebox*{!}{4cm}
	{\xymatrix{
	&F\times_B \widehat E\ar[rd]\ar[ld]\ar[rr]^\kappa&
	&E\times_B \widehat F\ar[rd]\ar[ld]&\\
	F\ar[rd]&&E\times_B \widehat E\ar[rd]\ar[ld]&&\widehat F\ar[ld]_{}\\
	&E\ar[rd]&&\widehat E\ar[ld]_{}&\\
	&&B&&
	}}.
	\end{eqnarray*}
By Proposition 
\ref{PropEqivalenceOfGerbes} 
the pair of the Pontrjagin duality triple has
the structure of a module pair over the
ring pair $F_{\widehat E^{\rm op}}\to B\times G\to B$, 
and  we can consider the Hilbert $C^*$-module 
$H(E,F)$ of this module pair over the $C^*$-algebra
$C^*(B\times G,F_{\widehat E^{\rm op}})$.

There is a second Hilbert $C^*$-module one is ought to consider.
Let $\Gamma_{0}(\widehat E,\widehat F^\C)$
be the bounded continuous  sections of the associated 
line bundle $\widehat F^\C$ of $\widehat F$
which vanish at fibre-wise infinity.
$\Gamma_{0}(\widehat E,\widehat F^\C)$ is a full Hilbert $C^*$-module
over the $C^*$-algebra $C_{0}(\widehat E)$ of continuous functions 
on $\widehat E$ which vanish at fibre-wise infinity.
The module structure is given by point-wise multiplication,
and the $C_{0}(\widehat E)$-valued inner product $\langle\ ,\ \rangle_{0}$
on $\Gamma_{0}(\widehat E,\widehat F^\C)$ is given 
for two sections $\hat\gamma,\hat \delta$ by
\begin{eqnarray}\label{EqThetaNichtTaeter}
	\langle\hat \gamma,\hat\delta\rangle_{0}:
	\widehat E\stackrel{\hat\gamma\times\hat\delta}{\longrightarrow}
	\widehat F^\C\times_{\widehat E} \widehat F^\C \stackrel{{\Theta}}{\longrightarrow} \C.
\end{eqnarray}
Here ${\Theta}:\widehat F^\C\times_{\widehat E} \widehat F^\C \to\C$ is
the canonical map $([\hat x,z],[\hat x\cdot t, z'] )\mapsto \overline{z} t z'$.
From a groupoid point of view $C_{0}(\widehat E)$ 
is the $C^*$-algebra of the $\U(1)$-central extension
\begin{eqnarray*}
\xymatrix{
\widehat E\!\times\!\U(1)\ \ar@<2pt>[d]\ar@<-2pt>[d]\ar@{>->}[r]
&\widehat E\times \U(1)\ar@{->>}[r]\ar@<2pt>[d]\ar@<-2pt>[d] 
&\widehat E\ar@<2pt>[d]\ar@<-2pt>[d]\\
\widehat E\ar[r]^=&\widehat E\ar[r]^=&\widehat E
},
\end{eqnarray*}
and $\Gamma_{0}(\widehat E,\widehat F^\C)$ is the 
Hilbert $C^*$-module of a $\U(1)$-equivariant Morita 
self-equivalence
 \begin{eqnarray*}
	\Big(
	\xymatrix{
	\widehat E\times\U(1) &
	\ar@{}[l]|{\hspace{0.7cm} \resizebox{!}{0,6cm}{\rotatebox {90}{{\rotatebox {90}{
	{\rotatebox {90}{$\circlearrowright$}}}}}}}
	F\ar[dl]\ar[dr] \ar@{}[r]|{\hspace{-0.5cm}\resizebox{!}{0,6cm}{\rotatebox {90}{$\circlearrowleft$}}}
	& \widehat E\times\U(1)\\
	\widehat E\ar@{<-}@<2pt>[u]\ar@{<-}@<-2pt>[u]&
	&\widehat E\ar@{<-}@<2pt>[u]\ar@{<-}@<-2pt>[u]
	}\Big)\ {\rotatebox {90}{$\circlearrowleft$}}\ \U(1).
\end{eqnarray*}

The following main theorem states 
that the Fourier transform identifies 
the two Hilbert $C^*$-modules we
are concerned with.

\begin{thm}\label{ThmTheMainThm}
	Let 
	\begin{eqnarray}\label{DiagPDTInTheMainThm}
	\resizebox*{!}{4cm}
	{\xymatrix{
	&F\times_B \widehat E\ar[rd]\ar[ld]\ar[rr]^\kappa&
	&E\times_B \widehat F\ar[rd]\ar[ld]&\\
	F\ar[rd]^q&&E\times_B \widehat E\ar[rd]\ar[ld]&&\widehat F\ar[ld]_{\hat q}\\
	&E\ar[rd]^p&&\widehat E\ar[ld]_{\hat p}&\\
	&&B&&
	}}
	\end{eqnarray}
	be a Pontrjagin duality triple, and let 
	\begin{eqnarray}\label{DiagBoehmerWald}
	\resizebox*{!}{4cm}
	{\xymatrix{
	&F_{\widehat E^{\rm op}}\times_B \widehat E\ar[rd]\ar[ld]\ar[rr]^{\kappa_{\widehat E}}&
	&\ar[rd]\ar[ld](B\times G)&\hspace{-1,6cm} \times_B \widehat E\times \U(1)\\
	F_{\widehat E^{\rm op}}\ar[rd]&&(B\times G)\times_B \widehat E\ar[rd]\ar[ld]&
	&\widehat E\times\U(1) \ar[ld]\\
	&B\times G\ar[rd]&&\widehat E\ar[ld]&\\
	&&B&&
	}}
	\end{eqnarray}
	be the Pontrjagin duality triple given by the bundle $\widehat E\to B$
	of (\ref{DiagPDTInTheMainThm}) as defined in (\ref{DiagTheFirstExaOfAPDT}).
	
	Then the Fourier transforms based on diagrams
	(\ref{DiagPDTInTheMainThm}) and (\ref{DiagBoehmerWald})
	extend to an isomorphism of Hilbert $C^*$-modules 
	\begin{eqnarray*}
		\Big( H(E,F),C^*(B\times G,F_{\widehat E^{\rm op}})\Big)
	&\cong& \Big(\Gamma_{0}(\widehat E,\widehat F^\C), C_{0}(\widehat E) \Big).\\
	(\gamma,\alpha)&\mapsto& (\hat\gamma,\hat \alpha)
	\end{eqnarray*}
	Moreover this isomorphism is natural in the base, i.e. it defines a functor 
	\begin{eqnarray*}
		{\rm Pon}^{\rm op}\to 
		({\rm Hilb}\text{-}C^*\text{-}{\rm mod})^{(0\cong1)},
	\end{eqnarray*}
	such that the value on the trivial Pontrjagin duality triple 
	over the one-point space is just the classical isomorphism of 
	Pontrjagin duality
	\begin{eqnarray*}
		(C^*(G),C^*(G))\cong (C_0(\widehat G),C_0(\widehat G)),
	\end{eqnarray*}
	where we regard the $C^*$-algebras as modules over themself. 
\end{thm}

\begin{pf}
	For two sections $\gamma,\delta:E\to F^\C$ let us compute 
	the Fourier transform of the inner product
	$\langle\gamma,\delta\rangle_c:B\times G\to F^\C_{\widehat E^{\rm op}}$.
	This is a function 
	$\widehat{\langle\gamma,\delta\rangle_c}:\widehat E\to\C$,
	namely
	\begin{eqnarray*}
		\widehat{\langle\gamma,\delta\rangle_c}(\hat e)
		&=&
		\int_G {\rm pr}_\C\big(\kappa^\C_{\widehat E}
		(\langle\gamma,\delta\rangle_c(b,g),	\hat e)\big)\ dg\\
		&=&
		\int_G {\rm pr}_\C\Big(\kappa^\C_{\widehat E}\big(
		\int_G \sigma_{\varrho}^\C(\gamma(e\cdot h),
		\delta(e\cdot(h+g))\ dh,\hat e\big)\Big)\ dg\\
		&=&
		\int_{G\times G} {\rm pr}_\C\Big(\kappa^\C_{\widehat E}\big(
		\sigma_{\varrho}^\C(\gamma(e\cdot h),
		\delta(e\cdot(h+g)),\hat e\big)\Big)\ d(g,h).
	\end{eqnarray*}
	To understand the integrated function we need to know the map 
	$\sigma_\varrho: F\times_{ E} F\to F_{\widehat E^{\rm op}}$
	explicitly. For $(x,x')\in  F\times_{ E} F$
	let $g\in G$ be such that $q(x')=q(x)\cdot g$,
	and for any $\hat e\in E|_{p(q(x))}$ let 
	$t\in\U(1)$ be such that
	${\rm pr}_{\widehat F}(\kappa(x',\hat e))=
	{\rm pr}_{\widehat F}(\kappa(x,\hat e))\cdot t$.
	Then, by use of (\ref{EqWerhettgeadcachtdat}), 
	we see that $\varrho(x,[\hat e,g,t])=x'$, i.e.  
	$\sigma_\varrho(x,x')=[\hat e,g,t]$.
	So for $[x,z]:=\gamma(e\cdot h)$ and
	$[x',z']:=\delta(e\cdot (h+g))$ we have
	$\sigma_\varrho^\C(\gamma(e\cdot h),\delta(e\cdot (h+g)))=
	[[\hat e,g,t],\overline z z']$. Thus the function we want to integrate 
	is 
	\begin{eqnarray*}
		&&
		{\rm pr}_\C\Big(\kappa^\C_{\widehat E}\big( \sigma_\varrho^\C
		(\gamma(e\cdot h),\delta(e\cdot (h+g))),\hat e\big)\Big)\\
		&=&{\rm pr}_\C\Big(\kappa^\C_{\widehat E}
		\big([[\hat e,g,t],\overline z z'],\hat e \big)\Big)\\
		&=&t \overline z z'\\
		&=&{\Theta} \big([{\rm pr}_{\widehat F}(\kappa(x,\hat e)),z],
		[{\rm pr}_{\widehat F}(\kappa(x',\hat e)),z']\big)\\
		&=&{\Theta} \Big({\rm pr}_{\widehat F^\C}
		\big(\kappa^\C(\gamma(e\cdot h),\hat e)\big),
		{\rm pr}_{\widehat F^\C}\big(\kappa^\C(\delta(e\cdot(h+g)',\hat e))\big)\Big)
	\end{eqnarray*}
	with $\Theta$ as in (\ref{EqThetaNichtTaeter}). 
	So if we do the $dg$-integration before the $dh$-integration in the above 
	integral, we see that 
	\begin{eqnarray*}
	\widehat{\langle\gamma,\delta\rangle_c}(\hat e)=
	\langle\hat \gamma,\hat\delta\rangle_0(\hat e).
	\end{eqnarray*}
	By similar arguments, we find that 
	the Fourier transform 
	$\hat{\ }:\Gamma_c(B\times G,F^\C_{\widehat E^{\rm op}})\to C_0(\widehat E)$
	is a morphism of algebras, i.e. $\widehat{\alpha\ast_\mu\beta}=\hat\alpha\hat \beta$, 
	and that the Fourier transform 
	$\hat{\ }:\Gamma_c(E,F^\C)\to \Gamma_0(\widehat E,\widehat F^\C)$
	is a morphism that preserves the  actions, i.e.
	$\widehat{\gamma\ast_\varrho\alpha}=\hat\gamma\hat\alpha$.
	By continuity, we get a map of Hilbert $C^*$-modules.

	The naturality of this map is obvious, and it 
	is also clear that over the one-point space we just 
	re\"obtain the classical Fourier transform which is 
	an isomorphism. 
	Using a partition of unity on the base space  $B$
	we find  that we have an isomorphism of Hilbert $C^*$-modules, 
	as it is an isomorphism locally. 
\end{pf}

\end{document}